%% file: NNCA.tex
\documentclass[manuscript,screen,review=false]{acmart}

\AtBeginDocument{%
  \providecommand\BibTeX{{%
    \normalfont B\kern-0.5em{\scshape i\kern-0.25em b}\kern-0.8em\TeX}}}

\setcopyright{acmcopyright}
\copyrightyear{2023}
\acmYear{2023}
\acmDOI{XXXXXXX.XXXXXXX}

\input{includes}

\begin{document}

\title{A new Nested Cross Approximation}

\author{Vaishnavi Gujjula}
\email{vaishnavihp@gmail.com}
\author{Sivaram Ambikasaran}
\email{sivaambi@iitm.ac.in}
\affiliation{%
  \institution{Indian Institute of Technology Madras}
  \city{Chennai}
  \state{Tamil Nadu}
  \country{India}
  \postcode{600036}
}


\renewcommand{\shortauthors}{Gujjula and Ambikasaran}

\begin{abstract}
In this article, we present a new Nested Cross Approximation (NNCA) for constructing $\mathcal{H}^{2}$ matrices. It differs from the existing NCAs~\cite{bebendorf2012constructing, zhao2019fast} in the technique of choosing pivots, a key part of the approximation. Our technique of choosing pivots is purely algebraic and involves only a single tree traversal. We demonstrate its applicability by developing a fast $\mathcal{H}^{2}$ matrix-vector product, that uses NNCA for the appropriate low-rank approximations. We illustrate the timing profiles and the accuracy of NNCA based $\mathcal{H}^{2}$ matrix-vector product. We also provide a comparison of NNCA based $\mathcal{H}^{2}$ matrix-vector product with the existing NCA based $\mathcal{H}^{2}$ matrix-vector products. A key observation is that NNCA performs better than the existing NCAs. In addition, using the NNCA based $\mathcal{H}^{2}$ matrix-vector product, we accelerate i) solving an integral equation in 3D and ii) Support Vector Machine (SVM). In the spirit of reproducible computational science, the implementation of the algorithm developed in this article is made available at \url{https://github.com/SAFRAN-LAB/NNCA}.
\end{abstract}



\maketitle

\section{Introduction}
In this article, we consider matrices $A\in\mathbb{R}^{N\times N}$ belonging to the class of $\mathcal{H}^{2}$
matrices~\cite{borm2003hierarchical,borm2003introduction,hackbusch2015hierarchical}.
Let the index sets of matrix $A$ be $I\times J$. For $i\in I$ and $j\in J$, let the $(i,j)^{th}$ entry of matrix $A$ be the evaluation of kernel function $K(x,y)$ at $x=p_{i}$ and $y=q_{j}$, where $\{p_{i}\}_{i\in I}$, and $\{q_{j}\}_{j\in J}$ are sets of points in $\mathbb{R}^{d}$, $K(x,y)\in \mathbb{R}$ and $K(x,y)$ is smooth everywhere except at $x=y$.
We denote the sets of points $\{p_{i}\}_{i\in I}$ and $\{q_{j}\}_{j\in J}$ by $P$ and $Q$ respectively.

Such matrices arise in many applications; Few of them are: particle simulations involving the Green's function of an elliptic PDE, discretization of an integral operator of an elliptic PDE, radial basis function interpolation, and covariance matrices in high dimensional statistics.

These matrices are usually large and dense, but certain sub-blocks of these matrices can be well-approximated by a low-rank matrix.
The construction of low-rank approximations of the appropriate matrix sub-blocks have been studied extensively in literature~\cite{cheng2005compression,dahmen2006compression,alpert1993wavelet,gimbutas2001coulomb,gimbutas2003generalized,martinsson2007accelerated,hackbusch2002h2,greengard1987fast,fong2009black,chan1987rank,zhao2005adaptive,rjasanow2002adaptive,ying2004kernel}.
One can classify the various low-rank approximations into two classes: analytic and algebraic.
Those methods that use certain analytic expansions of the underlying kernel function to construct low-rank approximations are termed analytic methods. And those methods that only need matrix entries and do not require knowledge of the underlying kernel are termed algebraic methods.
Examples of such analytic based low-rank construction include Taylor series expansions~\cite{greengard1987fast}, function interpolation~\cite{fong2009black}, etc.
Examples of algebraic based low-rank construction include rank revealing QR factorization~\cite{chan1987rank}, Adaptive Cross Approximation~\cite{zhao2005adaptive,rjasanow2002adaptive}, kernel-independent FMM~\cite{ying2004kernel}, etc.
There are a couple of advantages of algebraic methods over analytic methods: i) algebraic methods do not need analytic expansions of the kernel, instead only need matrix entries. ii) the ranks of the sub-blocks corresponding to far-field interactions obtained with algebraic methods are typically lower than that obtained with analytic methods, as algebraic methods are problem and domain-specific.

In this article, we develop an $\mathcal{O}(N)$, algebraic new Nested Cross Approximation (from now on abbreviated as NNCA). NCA, introduced in~\cite{bebendorf2012constructing}, an $\mathcal{O}(N\log(N))$ method, is a variant of Adaptive Cross Approximation~\cite{zhao2005adaptive,rjasanow2002adaptive} (from now on abbreviated as ACA) that provides \textit{nested bases} for $\mathcal{H}^{2}$ matrices~\cite{borm2003hierarchical,borm2003introduction,hackbusch2015hierarchical}.
The advantage of the nested bases in the context of $\mathcal{H}^{2}$ matrices is that it enables a linear time complexity for matrix-vector products (involving non-oscillatory kernels).

Our method differs from~\cite{bebendorf2012constructing} in two aspects:
\begin{enumerate}
  \item
  The technique of choosing self and far-field pivots: In~\cite{bebendorf2012constructing}, to find the far-field pivots, a geometric method is employed in a top-down fashion, which 
  chooses indices of points that are close to the tensor product Chebyshev nodes as pivots. And to find the self pivots both algebraic and geometric methods are described. 
  In this article, an algebraic method is employed in a bottom-up fashion, for both the self and far-field pivots, wherein the pivots of non-leaf cells at a parent level of the $2^{d}$ hierarchical tree are chosen from the pivots of cells at the child level. For more on self and far-field pivots, we refer the readers to Section~\ref{NCA}.
  \item
  The \textit{search space} of \textit{far-field pivots}: In~\cite{bebendorf2012constructing}, the search space of far-field pivots of a cluster of points is considered to be the entire \textit{far-field region} of the domain containing the support of the cluster of points. Whereas in this article the search space of far-field pivots of a cell belonging to the hierarchical tree is restricted to its interaction list (the terminology used in FMM).
\end{enumerate}

Zhao et al. in~\cite{zhao2019fast} developed an algebraic method of choosing pivots for NCA that has a time complexity of $\mathcal{O}(N)$.
The matrix partitioning and the search space of far-field pivots of~\cite{zhao2019fast} are similar to those presented in~\cite{bebendorf2012constructing}.
In~\cite{zhao2019fast}, a two-step process - a bottom-up followed by a top-down approach is employed to identify the pivots. 

NNCA, unlike the NCA presented in~\cite{zhao2019fast}, employs only a single tree traversal - a bottom-up traversal to identify the pivots.
The bottom-up traversal proposed in this article is similar to that of~\cite{zhao2019fast}, and we provide numerical evidence to illustrate that without compromising on accuracy one can eliminate the top-down traversal or the second step of~\cite{zhao2019fast}. 

The advantages of our method of choosing pivots over the ones in~\cite{bebendorf2012constructing,zhao2019fast} are:
\begin{enumerate}
\setlength\itemsep{0em}
  \item
  The search space of far-field pivots is smaller than those of the existing methods in~\cite{bebendorf2012constructing,zhao2019fast}.
  \item
  The nested bases 
  can be obtained from the pivot-choosing routine and do not need additional matrix entry evaluations for their computations.
\end{enumerate}
 As a consequence of these advantages, NNCA is computationally faster than the existing NCAs. In particular, the assembly time (the time taken to construct the $\mathcal{H}^{2}$ matrix representation) of NNCA is lower than that of the existing NCAs~\cite{bebendorf2012constructing,zhao2019fast}.


 We summarise the key aspects of this article here:
 \begin{enumerate}
\setlength\itemsep{0em}
   \item
   NNCA is proposed, where the pivot-choosing routine involves only a single tree traversal and the pivots are identified in a purely algebraic fashion.
   \item
   A comprehensive set of experiments to demonstrate the complexity and accuracy of the NNCA based $\mathcal{H}^{2}$ matrix-vector product are presented.
   \item
   A comparative study of various timing profiles and the accuracy of NNCA and NCAs is presented.
 \end{enumerate}

The rest of the article is organized as follows. In Section~\ref{section:lowRankRepresentation}, we present the new Nested Cross Approximation, wherein we detail the construction of $2^{d}$ tree upon which NNCA is built in Section~\ref{quadTree}, the admissibility condition for low-rank in Section~\ref{admissibilityCondition} and the construction of low-rank approximations in Section~\ref{NCA}. In Section~\ref{Bases}, we detail the construction of nested bases for NNCA. We describe the method to choose pivots for the NNCA in Section~\ref{Identification_of_Pivots}.
In Section~\ref{FMM_Matrix_Rep}, we present the steps to construct the $\mathcal{H}^{2}$ matrix representation.
In Section~\ref{sec:Algorithm}, we present the algorithm for NNCA based $\mathcal{H}^{2}$ matrix-vector product. We conclude the article with a comprehensive set of experiments in Section~\ref{NumericalResults} that provide various numerical benchmarks.

\section{New Nested Cross Approximation}
\label{section:lowRankRepresentation}

The key steps in constructing the new Nested Cross Approximation are: i) Sub-dividing the computational domain, or equivalently partitioning the matrix ii) Identifying the low-rank sub-blocks of matrix iii) Constructing low-rank approximations for the low-rank sub-blocks. We describe each of these steps in the sub-sections that follow.

\subsection{Construction of $2^{d}$ tree}\label{quadTree}
Let $D\in\mathbb{R}^{d}$ denote a compact hypercube in $d$ dimensions, amenable for constructing a $2^{d}$ hierarchical tree containing the support of the points in sets $P$ and $Q$.
To exploit the low-rank \textit{structure} of Hierarchical matrices, a hierarchical partitioning of the matrix $A$ into sub-blocks is needed. Equivalently, a hierarchical sub-division of the domain $D$ is to be performed.

A uniform $2^{d}$ tree is built in the domain $D$, upon which NNCA is built. Though an adaptive $2^{d}$ tree or a k-d tree too can be used, which partitions the matrix efficiently when the density of particles is widely varying, we consider a $2^{d}$ uniform tree for pedagogical reasons. It is to be noted that the algorithm described in this article is readily extendable to an adaptive tree or a k-d tree.
Level $0$ of the tree contains the \textit{root node} and is the domain $D$ itself. We sub-divide a cell (a node) at level $k$ into $2^{d}$ cells at level $k+1$. The former is said to be the parent of the latter and the latter are said to be the children of the former.
We stop further sub-division at level $\kappa$, if the number of points belonging to $P$ and the number of points belonging to $Q$ that lie in each cell of level $\kappa$ is less than or equal to $\nu$, where $\nu$ is a user-specified parameter that represents the maximum number of points in a leaf cell.
If we let $\mathcal{L}$ denote the set of all leaf cells, then
\begin{equation}
  D=\bigcup\{B:B\in \mathcal{L}\}.
\end{equation}

\subsection{Admissibility condition for low-rank}\label{admissibilityCondition}
For cells $X$ and $Y$ that belong to the $2^{d}$ tree, let $t^{X}$ and $s^{Y}$, defined below, denote the indices of points $\{p_{i}\}_{i\in I}$ and $\{q_{j}\}_{j\in J}$ that lie in cells $X$ and $Y$ respectively.
\begin{equation*}
  \begin{split}
    t^{X} &= \{i:i\in I \text{ and } p_{i}\in X\} \text{ and}\\
    s^{Y} &= \{j:j\in J \text{ and } q_{j}\in Y\}.
  \end{split}
\end{equation*}
We denote the matrix sub-block that captures the interaction between the clusters of points $\{p_{i}\}_{i\in t^{X}}$ and $\{q_{j}\}_{j\in s^{Y}}$ by $A_{t^{X}s^{Y}}$, whose $(i,j)^{th}$ entry is given by $A_{t^{X}s^{Y}}(i,j) = A(t^{X}(i),s^{Y}(j))$.

To construct the $\mathcal{H}^{2}$ matrix representation of $A$ using NNCA, we consider the following admissibility condition.
The matrix sub-block $A_{t^{X}s^{Y}}$ is approximated by a low-rank matrix if cells $X$ and $Y$ agree with the following admissibility condition for low-rank.
\begin{equation}
    max\{diam(X), diam(Y)\} \leq \eta dist(X,Y)
\end{equation}
where $\eta > 0$,
  \begin{align*}
    diam(X) &= sup\{\lVert x-y\rVert_{2}: x,y\in X\}, \\
    dist(X,Y) &= inf\{\lVert x-y\rVert_{2}: x\in X, y\in Y\}.
  \end{align*}
If cells $X$ and $Y$ satisfy the admissibility condition for low-rank, then the matrix sub-block $A_{t^{X}s^{Y}}$ is said to be \textit{admissible}.

\subsubsection{Preliminaries}
Before we proceed further we briefly explain ACA and present some notations 
that will be used in the rest of the article.

\noindent
\textbf{Adaptive Cross Approximation:}
The adaptive cross approximation (ACA) of an admissible matrix $A_{t^{X}s^{Y}}$, takes the form~\cite{bebendorf2000approximation,zhao2005adaptive}
\begin{equation}
A_{t^{X}s^{Y}}\approx A_{t^{X}s^{Y}}^{(k)} = UV^{\intercal} = A_{t^{X}\sigma^{Y}}A_{\tau^{X}\sigma^{Y}}^{-1}A_{\tau^{X}s^{Y}}
\end{equation}
where $\tau_{X}\subseteq t^{X}$ and $\sigma_{Y}\subseteq s^{Y}$ are termed the row and column pivots of the approximation. An iterative technique is used to identify the pivots in a heuristic fashion. In this article, we use the partially pivoted ACA algorithm~\cite{rjasanow2002adaptive}. For a given tolerance $\epsilon$, we stop further iterations when 
\begin{equation}
\|u_{k}\|_{2}\|v_{k}\|_{2}\leq \epsilon \|A_{t^{X}s^{Y}}^{(k)}\|_{F},
\end{equation}
where $u_{k}$ and $v_{k}$ are the $k^{th}$ column vectors of $U$ and $V$ matrices respectively and $A_{t^{X}s^{Y}}^{(k)}$ is the $k^{th}$ update of the approximation of $A_{t^{X}s^{Y}}$.

\noindent
\textbf{Notations:}

\begin{table}[h]
\centering
  \begin{tabular}{|l|p{0.8\textwidth}|}
    \hline
    $B$ & A cell (node) in the $2^{d}$-tree\\ \hline
    $\mathcal{C}(B)$ & $\{B':\text{ }B'\text{ is a child of }B\}$ \\ \hline
    $\mathcal{N}(B)$ & Set of neighbors of $B$, that consists of cells at the same tree level as $B$ which do not follow the admissibility condition for low-rank.\\ \hline
    $\mathcal{IL}(B)$ & Set of cells in the interaction list of cell $B$, that consists of children of $B$'s parent's neighbors that are not its neighbors.\\ \hline
    $\mathcal{A}(B)$ & Set of ancestors of $B$, that consists of all the nodes of the $2^{d}$-tree that lie on the path from root node to the node $B$.\\ \hline
  \end{tabular}
 \caption{Some notations that are followed in the rest of the article}
    \label{table:Notations}
 \end{table}

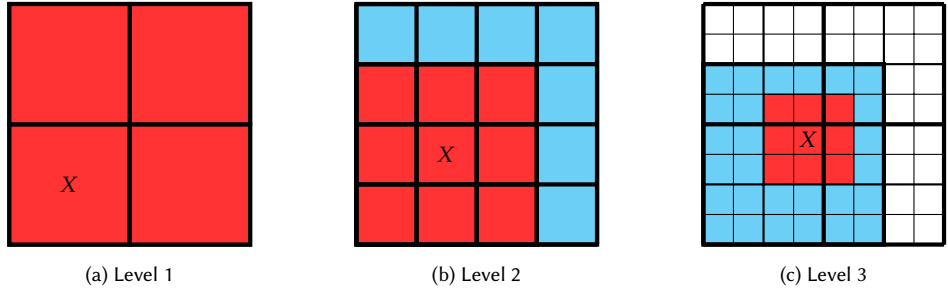
\begin{figure}[H]
\centering
	\begin{subfigure}[b]{0.3\textwidth}
\centering
\begin{tikzpicture}[scale=0.8]
    \draw[ultra thick] (-2,-2) rectangle (2,2);
    \draw[ultra thick,fill=red!80] (-2,-2) rectangle (0,0);
    \draw[ultra thick,fill=red!80] (0,-2) rectangle (2,0);
    \draw[ultra thick,fill=red!80] (0,0) rectangle (2,2);
    \draw[ultra thick,fill=red!80] (-2,0) rectangle (0,2);
    \node at (-1,-1) {$X$};
\end{tikzpicture}
\caption{Level 1}
	\end{subfigure}
	\begin{subfigure}[b]{0.3\textwidth}
\centering
\begin{tikzpicture}[scale=0.8]
    \draw[ultra thick,fill=cyan!50] (-2,-2) rectangle (2,2);
    \draw[ultra thick,fill=red!80] (-2,-2) rectangle (1,1);
    \draw[ultra thick,fill=cyan!50] (-2,-2) grid (2,2);
    \draw[ultra thick,fill=red!80] (-2,-2) grid (1,1);
    \node at (-0.5,-0.5) {$X$};
\end{tikzpicture}
\caption{Level 2}
    \end{subfigure}
    \begin{subfigure}[b]{0.3\textwidth}
\centering
\begin{tikzpicture}[scale=0.8]
\draw[ultra thick,fill=cyan!50] (-2,-2) rectangle (1,1);
    \draw[fill=red!80] (-1,-1) rectangle (0.5,0.5);
    \draw[step=2, ultra thick,fill=cyan!50] (-2,-2) grid (2,2);
    \draw[thick,fill=cyan!50] (-2,-2) grid (2,2);
    \draw[step=0.5, fill=red!80] (-2,-2) grid (2,2);
    \node at (-0.25,-0.25) {$X$};
\end{tikzpicture}
\caption{Level 3}
    \end{subfigure}
	\caption{Illustration of a cell $X$, its neighbors (the cells in red) and interaction list (the cells in cyan) at level 1, 2, and 3 of a quad-tree, with $\eta=\sqrt{2}$.}
	\label{fig:adm_eta_2}
\end{figure}

We illustrate in Figure~\ref{fig:adm_eta_2}, the neighbors and interaction list of a cell with $\eta=\sqrt{2}$.

We term the \textit{far-field region} of a cell $X$, represented by $F(X)$, as
\begin{equation*}
  F(X) = \bigcup_{X'\in\mathcal{A}(X)} \mathcal{IL}(X')
\end{equation*}

\subsection{Construction of low-rank approximations}\label{NCA}
Consider two cells $X$ and $Y$, such that $Y\in \mathcal{IL}(X)$.
The admissible block $A_{t^{X}s^{Y}}$ is approximated by a low-rank matrix as
\begin{equation} \label{eq:NCA}
  A_{t^{X}s^{Y}} \approx \underbrace{A_{t^{X}s^{X,i}} (A_{t^{X,i}s^{X,i}})^{-1}}_{U_{X}} \underbrace{A_{t^{X,i}s^{Y,o}}}_{S_{X,Y}} \underbrace{(A_{t^{Y,o}s^{Y,o}})^{-1} A_{t^{Y,o}s^{Y}}}_{V_{Y}^{\intercal}}
\end{equation}
where $t^{X,i}\subset t^{X}$, $s^{X,i}\subset \mathcal{F}^{X,i}$, $t^{Y,o}\subset \mathcal{F}^{Y,o}$ and $s^{Y,o}\subset s^{Y}$.
$\mathcal{F}^{X,i}$ and $\mathcal{F}^{Y,o}$ are defined as
\begin{equation}\label{eq:Fj}
  \mathcal{F}^{X,i} = \{s^{X'}:X'\in \mathcal{IL}(X)\} \text{ and}
\end{equation}
\begin{equation}\label{eq:Fi}
  \mathcal{F}^{Y,o} = \{t^{Y'}:Y'\in \mathcal{IL}(Y)\},
\end{equation}
whereas in~\cite{bebendorf2012constructing,zhao2019fast} they are defined as
\begin{equation}\label{eq:Fj_old}
  \mathcal{F}^{X,i} = \{s^{X'}:X'\in F(X)\} \text{ and}
\end{equation}
\begin{equation}\label{eq:Fi_old}
  \mathcal{F}^{Y,o} = \{t^{Y'}:Y'\in F(Y)\}.
\end{equation}

We term $t^{X,i}$ and $s^{X,i}$ as the \textit{incoming row pivots} and \textit{incoming column pivots} of $X$ respectively.
$t^{Y,o}$ and $s^{Y,o}$ are termed the \textit{outgoing row pivots} and \textit{outgoing column pivots} of $Y$ respectively.
Pivots $t^{X,i}$ and $s^{Y,o}$ represent $t^{X}$ and $s^{Y}$ respectively.
Pivots $s^{X,i}$ and $t^{Y,o}$ represent $\mathcal{F}^{X,i}$ and $\mathcal{F}^{Y,o}$ respectively.
Though the points indexed by $s^{X,i}$ do not lie in $X$, $s^{X,i}$ is considered to depend only on $X$ as it represents the points (belonging to $Q$) lying in the interaction-list region of $X$. Similarly $t^{Y,o}$ is considered to depend only on $Y$.

$U_{X}$ and $V_{Y}^{\intercal}$, pointed in equation~\eqref{eq:NCA}, are termed the \textit{column basis} of $X$ and \textit{row basis} of $Y$ respectively.

If $A_{t^{X}s^{Y}}$ is admissible, it follows that $A_{t^{Y}s^{X}}$ is also admissible. The low rank approximation of $A_{t^{Y}s^{X}}$ via NNCA takes the form
\begin{equation} \label{eq:NCA2}
  A_{t^{Y}s^{X}} \approx \underbrace{A_{t^{Y}s^{Y,i}} (A_{t^{Y,i}s^{Y,i}})^{-1}}_{U_{Y}} \underbrace{A_{t^{Y,i}s^{X,o}}}_{S_{Y,X}} \underbrace{(A_{t^{X,o}s^{X,o}})^{-1} A_{t^{X,o}s^{X}}}_{V_{X}^{\intercal}}
\end{equation}
where $t^{Y,i}\subset t^{Y}$, $s^{Y,i}\subset \mathcal{F}^{Y,i}$, $t^{X,o}\subset \mathcal{F}^{X,o}$ and $s^{X,o}\subset s^{X}$.

Therefore a cell $X$ is associated with four sets of pivots $t^{X,i}$, $s^{X,i}$, $t^{X,o}$ and $s^{X,o}$ and bases $U_{X}$ and $V_{X}^{\intercal}$.
$s^{X,i}$ and $t^{X,o}$ represent the points (belonging to $Q$ and $P$ respectively) that lie in the far-field region of $X$ and hence are also termed the \textit{far-field pivots} of $X$.
$t^{X,i}$ and $s^{X,o}$ represent the points (belonging to $P$ and $Q$ respectively) that lie in $X$ and hence are also termed the \textit{self pivots} of $X$. For details on error estimates of equation~\eqref{eq:NCA} we refer the readers to~\cite{bebendorf2012constructing}.

With reference to equation~\eqref{eq:NCA}, the contribution of particles indexed by $s^{Y}$ at $t^{X}$ is captured via the contribution of particles indexed by $s^{Y,o}$ at $t^{X,i}$.
Here $U_{X}$ can be interpreted as an algebraic interpolation operator, interpolating from $t^{X,i}$ to $t^{X}$.
Similarly $V_{Y}^{\intercal}$ can be interpreted as an algebraic anterpolation operator, anterpolating from $s^{Y}$ to $s^{Y,o}$.

An interpolation operator can therefore be constructed which interpolates from the incoming row pivots of a non-leaf cell to the incoming row pivots of its children.
Similarly, an anterpolation operator can be constructed which anterpolates from the outgoing column pivots of child cells to the outgoing column pivots of their parent.
This gives a technique to construct nested bases.

\subsubsection{Construction of Nested Bases}\label{Bases}
To construct the low-rank approximation of the far-field interactions of $B$, one needs matrices $U_{B}$ and $V_{B}^{\intercal}$.

For a leaf cell $B$,
    \begin{equation}
    U_{B} := A_{t^{B}s^{B,i}} (A_{t^{B,i} s^{B,i}})^{-1} \hspace{5mm} \text{and} \hspace{5mm}
    V_{B}^{\intercal} := (A_{t^{B,o} s^{B,o}})^{-1} A_{t^{B,o}s^{B}}.
    \end{equation}
Matrices $U_{B}$ and $V_{B}^{\intercal}$, the column and row bases of $B$ respectively, are also termed the L2P (local-to-particle) and P2M (particle-to-multipole) translation operators of $B$ (the terminology used in FMM).

The column and row bases of non-leaf cells are constructed in a nested fashion: The bases of a cell are expressed in terms of the bases of its children.
For a non-leaf cell $B$,
  \begin{equation}
    U_{B} = \begin{bmatrix}
            U_{B_{1}} & 0 & \hdots & 0\\
            0 & U_{B_{2}} &  & 0\\
            \vdots &  & \ddots & \\
            0 & 0 & \hdots & U_{B_{2^{d}}}
            \end{bmatrix}
            \begin{bmatrix}
            C_{B_{1}B}\\
            C_{B_{2}B}\\
            \vdots\\
            C_{B_{2^{d}}B}\\
            \end{bmatrix},
  \end{equation}

            \begin{equation}
            V_{B} = \begin{bmatrix}
                    V_{B_{1}} & 0 & \hdots & 0\\
                    0 & V_{B_{2}} &  & 0\\
                    \vdots &  & \ddots & \\
                    0 & 0 & \hdots & V_{B_{2^{d}}}
                    \end{bmatrix}
                    \begin{bmatrix}
                    T_{BB_{1}}^{\intercal}\\
                    T_{BB_{2}}^{\intercal}\\
                    \vdots\\
                    T_{BB_{2^{d}}}^{\intercal}\\
                    \end{bmatrix},
  \end{equation}
 where $\{B_{c}\}_{c=1}^{2^{d}}\in \mathcal{C}(B)$ and
matrices $\{C_{B_{c}B}\}_{1\leq c \leq 2^{d}}$ and $\{T_{BB_{c}}^{\intercal}\}_{1\leq c \leq 2^{d}}$ take the following form
  \begin{equation}\label{eq:S_Non_Directional_L2L}
    C_{B_{c}B} =  A_{t^{B_{c},i}s^{B,i}} (A_{t^{B,i} s^{B,i}})^{-1} \hspace{1mm} \text{and} \hspace{1mm}
    T_{BB_{c}} = (A_{t^{B,o} s^{B,o}})^{-1} A_{t^{B,o}s^{B_{c},o}} ,
  \end{equation}
$\forall c\in\{1,2,...,2^{d}\}$. Matrices $\{C_{B_{c}B}\}_{1\leq c \leq 2^{d}}$ and $\{T_{BB_{c}}^{\intercal}\}_{1\leq c \leq 2^{d}}$ are termed the column translation matrices or the L2L's (local-to-local) and row translation matrices or M2M's (multipole-to-multipole) of $B$ respectively (the terminology used in FMM).

For a cell $B$, $U_{B}$ and $V_{B}^{\intercal}$ approximate $A_{t^{B}s^{B,i}} (A_{t^{B,i} s^{B,i}})^{-1}$ and \\$(A_{t^{B,o} s^{B,o}})^{-1} A_{t^{B,o}s^{B}}$.  We refer the readers to~\cite{bebendorf2012constructing} for the error estimates.
\subsubsection{Identification of Pivots}\label{Identification_of_Pivots}
To obtain the low-rank approximations of the admissible sub-blocks of the matrix, it remains to discuss how the pivots are selected.
For this, one needs four sets of pivots $t^{B,i}$, $s^{B,i}$, $t^{B,o}$, and $s^{B,o}$, defined for each cell $B$ of the $2^{d}$ tree.

$t^{B,i}$ and $s^{B,o}$ are chosen from $t^{B}$ and $s^{B}$ respectively, and hence the \textit{search space} of self pivots of a cell is itself.

Bebendorf et al. in~\cite{bebendorf2012constructing} considered the search space of far-field pivots of a cluster of points to be the entire far-field region of the domain containing the support of the cluster of points. Zhao et al. in~\cite{zhao2019fast} follow a two-stage process to find pivots. Stage 1 computes the local pivots from the local far-field region (or the interaction-list region) and stage 2 uses the local pivots as input and finds the pivots corresponding to the entire far-field region. Hence Zhao et al. too consider the search space of far-field pivots of a cluster of points to be the entire far-field region of the domain containing the support of the cluster of points.

In this article, the far-field pivots, $s^{B,i}$ and $t^{B,o}$, are chosen from $\mathcal{F}^{Bi}$ and $\mathcal{F}^{Bo}$ respectively, which contain indices of points (belonging to $P$  and $Q$ respectively) that lie in the interaction list of $B$. So the \textit{search space} of far-field pivots of a cell is its interaction list region.
We claim based on numerical evidence (refer Section~\ref{NumericalResults}) that it is sufficient to choose the far-field pivots of a cell from the indices of points in its interaction list region. 

In Subsections~\ref{ComparisonExperiment} and~\ref{sec:exp2} we illustrate the convergence of NNCA. We also compare the NNCA with that of the existing NCAs in Subsection~\ref{ComparisonExperiment}, and it is to be observed that NNCA performs better than the existing NCAs. 

We now describe the method to choose pivots of all cells in the $2^{d}$ tree in a nested fashion, where we obtain pivots of cells at a parent level from the pivots of cells at the child level.
We traverse up the tree (starting at the leaf level) in a reverse level-order fashion to find pivots of all cells using the two steps given below.
 \begin{enumerate}
   \item
    Construct sets $\tilde{t}^{B,i}$, $\tilde{s}^{B,i}$, $\tilde{t}^{B,o}$ and $\tilde{s}^{B,o}$ that represent $t^{B}$, $\mathcal{F}^{B,i}$, $\mathcal{F}^{B,o}$ and $s^{B}$ respectively.\\
     For a leaf cell $B$, construct sets
    \begin{equation}
      \tilde{t}^{B,i}:=t^{B}, \hspace{5mm}  \hspace{5mm} \tilde{s}^{B,i}:=\bigcup_{B'\in \mathcal{IL}(B)} s^{B'},
    \end{equation}
    \begin{equation}
       \tilde{t}^{B,o}:=\bigcup_{B'\in \mathcal{IL}(B)}t^{B'} \hspace{5mm} \text{and} \hspace{5mm}
       \tilde{s}^{B,o}:=s^{B}.
    \end{equation}

    For a non-leaf cell $B$, construct sets
    \begin{equation}
      \tilde{t}^{B,i}:=\bigcup_{B'\in\mathcal{C}(B)} t^{B',i}, \hspace{5mm}  \hspace{5mm} \tilde{s}^{B,i}:=\bigcup_{B'\in \mathcal{IL}(B)} \bigcup_{B''\in \mathcal{C}(B')} s^{B'',o},
    \end{equation}
    \begin{equation}
       \tilde{t}^{B,o}:=\bigcup_{B'\in \mathcal{IL}(B)} \bigcup_{B''\in \mathcal{C}(B')} t^{B'',i} \hspace{5mm} \text{and} \hspace{5mm}
       \tilde{s}^{B,o}:=\bigcup_{B'\in\mathcal{C}(B)}s^{B',o}.
    \end{equation}

\item
 Perform ACA~\cite{zhao2005adaptive,rjasanow2002adaptive} on the matrix $A_{\tilde{t}^{Bi}\tilde{s}^{Bi}}$ with tolerance $\epsilon_{NCA}$.
 The row and column pivots chosen by ACA are then assigned to pivots $t^{Bi}$ and $s^{Bi}$ respectively.
Similarly, perform ACA on the matrix $A_{\tilde{t}^{Bo}\tilde{s}^{Bo}}$ to get the pivots $t^{Bo}$ and $s^{Bo}$.

 \end{enumerate}

\begin{remark}
In addition to the difference in search space of far-field pivots, our technique of choosing pivots differs from the ones in~\cite{bebendorf2012constructing, zhao2019fast} in
the method of construction of far-field pivots. We follow a bottom-up approach: We start at the leaf level and traverse up the tree to find the far-field pivots in a recursive manner, wherein the far-field pivots at a parent level are constructed from the far-field pivots at the child level as explained above.
In~\cite{bebendorf2012constructing}, a top-down approach is followed: The far-field pivots at a child level are constructed from the far-field pivots at the parent level.
In~\cite{zhao2019fast}, a bottom-up approach followed by a top-down approach is employed, wherein the bottom-up approach is similar to the one proposed in this article. For a cell $B$,
the bottom-up approach is used to choose the partial far-field pivots (or local far-field pivots) from its interaction list region\footnote{The matrix partitioning in~\cite{zhao2019fast} is different from that presented in this article. But we use the notion of cell, though is incorrect, to present the technicalities without introducing additional notations.}. The top-down approach is used to find additional far-field pivots from the interaction list of its ancestors other than itself. 
In Section~\ref{NumericalResults}, we demonstrate numerically that the top-down approach of~\cite{zhao2019fast} can be avoided with no compromise on accuracy.
 \end{remark}
\begin{remark}
There are a couple of advantages of our method over the existing methods:

i) The search space of far-field pivots is smaller than those considered in~\cite{bebendorf2012constructing,zhao2019fast}. Smaller search space leads to applying ACA on smaller matrix sizes, in the second step of the method. As a result, our method is computationally faster than the existing methods without any substantial compromise in accuracy.
The numerical results we provide in Section~\ref{ComparisonExperiment} demonstrate the timing profiles and accuracy of our method.

ii)
For a cell $B$, the second step of the pivot-choosing routine involves using ACA which takes the form
\begin{equation*}
  A_{\tilde{t}^{B,i}\tilde{s}^{B,i}}\approx A_{\tilde{t}^{B,i}s^{B,i}}(A_{t^{B,i}s^{B,i}})^{-1}A_{t^{B,i}\tilde{s}^{B,i}}.
\end{equation*}
The advantage of using ACA in NNCA and existing NCAs is that the LU decomposition of $(A_{t^{B,i}s^{B,i}})^{-1}$, is available as a byproduct of ACA~\cite{bebendorf2009recompression}.
In addition to this advantage, NNCA benefits from ACA in the evaluation of $U_{B}$ for leaf cells $B$, and $C_{B_{c}B}$ matrices for non-leaf cells $B$, where $\{B_{c}\}_{c=1}^{2^{d}}\in \mathcal{C}(B)$.
For a non-leaf cell $B$, $A_{\tilde{t}^{B,i}s^{B,i}}(A_{t^{B,i}s^{B,i}})^{-1}\equiv \begin{bmatrix}
C_{B_{1}B}^{\intercal} & C_{B_{2}B}^{\intercal} & \hdots & C_{B_{2^{d}}B}^{\intercal}
\end{bmatrix}^{\intercal}
$.
For a leaf cell $B$, $A_{\tilde{t}^{B,i}s^{B,i}}(A_{t^{B,i}s^{B,i}})^{-1}\equiv  U_{B}$.
Hence ACA enables us to evaluate $U_{B}$ for leaf cells $B$, and $C_{B_{c}B}$ matrices for non-leaf cells $B$,
with no additional matrix entry evaluations. This does not hold true for the existing methods.
A similar observation is to be made with the $V_{B}$ and $T_{BB_{c}}$ matrices.
\end{remark}

\section{$\mathcal{H}^{2}$ matrix representation}\label{FMM_Matrix_Rep}
The construction of $\mathcal{H}^{2}$ matrix representation involves obtaining low-rank approximations of matrix sub-blocks $A_{t^{X}s^{Y}}$, for all cells $X$ at all levels of the $2^{d}$ tree, where $Y\in\mathcal{IL}(X)$. And the rest of the matrix sub-blocks are built exactly (up to roundoff) with no compression involved.
The algorithm to construct the $\mathcal{H}^{2}$ matrix representation is described below.
\begin{enumerate}
  \item
  Compute pivots and the column basis or column translation matrices and the row basis or row translation matrices of all cells $X$ at all levels of the tree as described in Section~\ref{Identification_of_Pivots}.
  \item
  Construct matrices $S_{X,Y}=A_{t^{Xi}s^{Yo}}$, for all cells $X$ at all levels of the $2^{d}$ tree where $Y\in\mathcal{IL}(X)$.
  \item
  Construct matrices $A_{t^{X}s^{Y}}$ for all leaf cells $X$ of the $2^{d}$ tree where $Y\in\mathcal{N}(X)$.
\end{enumerate}

\section{$\mathcal{H}^{2}$ matrix-vector product}
\label{sec:Algorithm}
In this section, NNCA based $\mathcal{H}^{2}$ matrix-vector product to evaluate $u$,
\begin{equation}
  u=Aw,
\end{equation}
is presented, where $u,w\in\mathbb{R}^{N}$ and $A\in\mathbb{R}^{N\times N}$.

For a cell $B$, let $w^{B}$ and $u^{B}$ be defined as
\begin{equation*}
    w^{B} = [w_{j_{1}}, w_{j_{2}},\text{\textellipsis},w_{j_{b_{1}}}] \text{ where } \{j_{c}\}_{1\leq c\leq b_{1}}= s^{B}, \text{ } b_{1}=\lvert s^{B}\rvert\text{ and}\\
\end{equation*}
\begin{equation*}
    u^{B} = [u_{i_{1}}, u_{i_{2}},\text{\textellipsis},u_{i_{b_{2}}}] \text{ where } \{i_{c}\}_{1\leq c\leq b_{2}}= t^{B}, \text{ }b_{2}=\lvert t^{B}\rvert.
\end{equation*}

The algorithm is as follows:
\begin{enumerate}
  \item
  \textbf{NNCA: }Construct the $\mathcal{H}^{2}$ matrix representation of $A$, as described in Section~\ref{FMM_Matrix_Rep}.
  \item
  \textbf{Upward Pass: }For all leaf cells $B$, compute 
  \begin{equation*}
    w^{B,o} = V^{\intercal}_{B}w^{B}.
  \end{equation*}
  For all non-leaf cells $B$ at level $k$, compute 
  \begin{equation*}
    w^{B,o} = \sum_{B'\in \mathcal{C}(B)} T^{\intercal}_{BB'}w^{B',o}.
  \end{equation*}
by recursion, $\kappa-1\geq k\geq 0$
  \item
  \textbf{Transverse Pass:} For all cells $B$ at all levels, compute 
  \begin{equation*}
    u^{B,i} = \sum_{B'\in\mathcal{IL}(B)}S_{BB'} w^{B',o}
  \end{equation*}
  \item
  \textbf{Downward Pass: }For all non-leaf cells $B'$ at level $k$, compute
  \begin{equation*}
    u^{B',i} := u^{B',i} + C_{B'B} u^{B,i}
  \end{equation*}
 by recursion, $1\leq k\leq \kappa$,
 where $B$ is parent of $B'$.

  For all leaf cells $B$, compute
  \begin{equation*}
    u^{B} := U_{B} u^{B,i}.
  \end{equation*}

  \item
  For all leaf cells $B$, add the near field interaction to $u^{B}$
  \begin{equation*}
    u^{B} := u^{B} + \sum_{B'\in\mathcal{N}(B)}A_{t^{B}s^{B'}}w_{B'}.
  \end{equation*}
\end{enumerate}

\subsection{Time complexity}
\noindent
  \textbf{NNCA.}  
\noindent
   Let $a=max\{\bigcup\{\{\lvert t^{Bi}\rvert, \lvert t^{Bo}\rvert, \lvert s^{Bi}\rvert, \lvert s^{Bo}\rvert\}:B\text{ belongs to quad-tree}\}\}$. Assume the maximum leaf size (maximum number of particles a leaf can have) $\nu=\mathcal{O}(a)$. Then for a cell $B$ in the $2^{d}$ tree:
   $\lvert \tilde{t}^{Bi}\rvert=\mathcal{O}(a)$, $\lvert \tilde{t}^{Bo}\rvert=\mathcal{O}(a)$, $\lvert \tilde{s}^{Bi}\rvert=\mathcal{O}(a)$ and $\lvert\tilde{s}^{Bo}\rvert=\mathcal{O}(a)$.
   \begin{itemize}
       \item 
   The pivots of a cell are chosen from either the pivots of its children (for a non-leaf cell), of cardinality $n_{1} = \mathcal{O}(a)$, or the particles lying in the cell (for a leaf cell), of cardinality $n_{2} = \nu = \mathcal{O}(a)$. So, the complexity of finding pivots of cell $B$ is $\mathcal{O}(a^{3})$, since performing ACA on a matrix of size $n_{1}\times n_{2}$ and numerical rank\footnote{Numerical rank of a matrix $A$, $r_{\epsilon}(A)$, is defined as $\min\{k\in\{1,2,...,N\}:\frac{\sigma_{k}}{\sigma_{1}}<\epsilon\}$, where $\sigma_{1}, \sigma_{2}, ..., \sigma_{N} \text{ are the singular values of matrix } A \text{ and } \sigma_{1}\geq\sigma_{2}\geq...\sigma_{N}$} $\mathcal{O}(a)$, costs $(n_{1}+n_{2})a^{2}$. 
       \item 
   It is to be noted that no additional cost is needed to compute the column basis or column translation matrices and row basis or row translation matrices of $B$, as they can be obtained as byproducts of the pivot-choosing routine.
       \item 
   Assuming a matrix entry can be obtained in $\mathcal{O}(1)$ time, the complexity of computing matrices $\{S_{B,X}:X\in \mathcal{IL}(B)\}$ of cell $B$ is $\mathcal{O}(a^{2})$. 
       \item 
   The complexity of obtaining the matrices $\{A_{t^{B}s^{B'}}: B'\in\mathcal{N}(B)\}$ is $\mathcal{O}(a^{2})$.
\end{itemize}
Hence the sum of the complexities of finding pivots and computing the column basis or column translation matrices, row basis or row translation matrices, matrices $\{S_{B,X}:X\in \mathcal{IL}(B)\}$, matrices $\{A_{t^{B}s^{B'}}: B'\in\mathcal{N}(B)\}$ is $\mathcal{O}(a^{3})$. For a tree with maximum leaf size $\nu$, the number of cells at all levels is equal to $\mathcal{O}(N/\nu) = \mathcal{O}(N/a)$.
   So the total complexity of finding pivots and computing the column basis or column translation matrices, row basis or row translation matrices, matrices $\{S_{B,X}:X\in \mathcal{IL}(B)\}$, matrices $\{A_{t^{B}s^{B'}}: B'\in\mathcal{N}(B)\}$ of all cells $B$ at all levels of the $2^{d}$ tree or equivalently, the complexity of computing the $\mathcal{H}^{2}$ representation of the matrix is $\mathcal{O}(a^{2}N)$.
   
   The kernel functions considered in this article are asymptotically smooth away from the singularity, ($\lVert x-y\rVert_{2}=0$), and hence the matrices corresponding to the far-field interactions can be efficiently approximated by a low-rank matrix, whose rank is independent of $N$~\cite{hackbusch2002blended,brandt1998multilevel}. As a result, the upper bound of the rank of the far-field interaction, $a$, is independent of $N$. Hence the time complexity of constructing NNCA is $\mathcal{O}(N)$.

  \bigskip
  \noindent
  \textbf{$\mathcal{H}^{2}$ matrix-vector product.}
  Step $1$ of the algorithm, the construction of NNCA, costs $\mathcal{O}(N)$, as stated above. 
  Steps $2-5$ of the algorithm are the usual steps in a $\mathcal{H}^{2}$ matrix-vector product algorithm, which cost $\mathcal{O}(N)$.
Hence the overall time complexity of the $\mathcal{H}^{2}$ matrix-vector product algorithm is $\mathcal{O}(N)$.

\subsection{Memory complexity}
\bigskip
\noindent
  \textbf{NNCA.}
For a cell $B$ in the $2^{d}$ tree, the cost of storing the column basis /column translation matrices, matrices $\{S_{B,X}:X\in \mathcal{IL}(B)\}$, and row basis/ row translation matrices is $\mathcal{O}(a^{2})$. The cost of storing these matrices for all cells of the $2^{d}$ tree is $\mathcal{O}(aN)$. On similar lines as that of the total time complexity, the total memory complexity simplifies to $\mathcal{O}(N)$. 
   
\section{Numerical Results}\label{NumericalResults}
We perform a total of six experiments to demonstrate the performance of NNCA in 2D, 3D and 4D:
\begin{enumerate}
    \item
     Matrix-vector product with uniform distribution of particles in 2D and its comparison with the existing NCAs.
    \item
     Matrix-vector product with a non-uniform distribution of points in 2D.
    \item
     Matrix-vector product with a uniform distribution of points in 3D.
    \item
     Matrix-vector product with a non-uniform distribution of points in 3D.
    \item
     Integral equation solver in 3D.
    \item
     Kernel SVM (Support Vector Machine) in 2D and 4D.
\end{enumerate}
 
In all the experiments we use $\eta=\sqrt{2}$.
In experiments 1-4, we consider kernel functions of the form: $K(x,y) = K(r)$ $(r=\lVert x-y\rVert_{2})$, where the vector to be applied to the matrix, $w$, is taken to be a random vector. 

In all the experiments, $P$ and $Q$ are considered to be the same, and the kernels dealt with are symmetric. So the matrix $A$ is symmetric. As a result its sufficient to compute pivots $t^{Bi}$ and $s^{Bi}$ and the assignments $t^{Bo}:=s^{Bi}$ and $s^{Bo}:=t^{Bi}$ follow as the premise. Also, for non-leaf cells $B$, its sufficient to compute the $C_{B'B}$ operator $\forall B'\in\mathcal{C}(B)$. Similarly for leaf cells $B$ it is sufficient to compute $U_{B}$. The assignments $T_{BB'}:=C_{B'B}$ for non-leaf cells and $V_{B} := U_{B}$ for leaf cells follow as the premise.

In experiment 1, for the purpose of comparison with the existing methods, we implemented the NCAs of~\cite{bebendorf2012constructing} and~\cite{zhao2019fast} as well.
For the implementation of the NCA of~\cite{bebendorf2012constructing}, the cardinality of $\tilde{t}^{Bi}$, for all cells $B$ at all levels, is set to
$min\{k_{\epsilon}^{2},\lvert t^{B}\rvert\}$,
where $k_{\epsilon}$ denotes the number of terms in the truncated Taylor series expansion of the kernel function with a relative error bound of $O(\epsilon)$.
For a $\mathcal{H}^{2}$ matrix with $\eta=\sqrt{2}$,
$k_{\epsilon}$ is $\lceil(-\log_{c}(\epsilon))\rceil$, where $c\approx 1.828$~\cite{greengard1988rapid}. 
In this article we choose $k_{\epsilon}=\lceil(-\log_{1.25}(\epsilon))\rceil > \lceil(-\log_{c}(\epsilon))\rceil$.
It is to be observed from Figures~\ref{1VsN} and~\ref{2VsN} that at large values of $N$, even with a high value of $k_{\epsilon}$
the relative error obtained with Bebendorf et al.'s algorithm is larger than that of NNCA and Zhao et al.'s NCA. The notations described in Table~\ref{table:ResultsNotations} will be used in the rest of the section.

\begin{table}[h]
  \centering
  \begin{tabular}{|l|p{0.8\textwidth}|}
    \hline
 $N$ & System size \\
 \hline
 mem. & Memory needed in GB to store the matrix in NCA/NNCA based $\mathcal{H}^{2}$ matrix representation\\
 \hline
 T\textsubscript{a} & Assembly time, that is the time taken in seconds to construct the NCA/NNCA based $\mathcal{H}^{2}$ matrix representation\\
 \hline
 T\textsubscript{m} & Time taken in seconds to compute matrix-vector product using NCA/NNCA\\
 \hline
 $\epsilon_{m}$ & Relative error in the matrix-vector product in 2-norm sense\\
 \hline
 T\textsubscript{s} & Time taken in seconds to solve using an iterative solver\\
 \hline
 $\epsilon_{s}$ & Relative error in the solution in 2-norm sense\\
 \hline
  iter. & Number of iterations taken by an iterative solver to converge to a given accuracy.\\
 \hline
\end{tabular}
\caption{List of notations followed in the rest of the section}
   \label{table:ResultsNotations}
\end{table}

The algorithm is implemented in C++. Experiments 1 and 2 were run on a 2.3GHz Intel Core i5 processor with 4 Openmp threads. Experiments 3, 4, and 5 were run on an Intel Xeon Gold, 2.5 GHz processor with 8 OpenMP threads. Experiment 6 is run on a 2.3GHz Intel Core i5 processor with no parallelization.

In the spirit of reproducible computational science, the implementation of the algorithm developed in this article is made available at \url{https://github.com/SAFRAN-LAB/NNCA}. The documentation of this library together with the data and inputs that reproduce the results illustrated in this section is available at \url{https://nnca.readthedocs.io/en/latest/}.

\subsection{Experiment 1: Matrix-vector product with uniform distribution of particles in 2D and its comparison with the existing NCAs}\label{ComparisonExperiment}
In this sub-section, we present various benchmarks for NCA and NNCA based $\mathcal{H}^{2}$ matrix-vector product in 2D. 
mem., $T_{a}$, $T_{m}$, and $\epsilon_{m}$ are compared with those of the existing methods by Bebendorf et al.~\cite{bebendorf2012constructing} and Zhao et al.~\cite{zhao2019fast}. 
We experiment with two kernels: i) $\left(\frac{r(\log(r)-1)}{a(\log(a)-1)}\right)\chi_{r<a}+\left(\frac{\log(r)}{\log(a)}\right)\chi_{r\geq a}$ ii) $\left(\frac{r}{a}\right)\chi_{r<a}+\left(\frac{a}{r}\right)\chi_{r\geq a}$,
where $a$ is set to $10^{-4}$ and the particles are distributed uniformly in the domain $[-1,1]^{2}$.

\begin{figure}[H]
  \begin{center}

      \begin{subfigure}[b]{0.45\textwidth}
        \centering
    \includegraphics{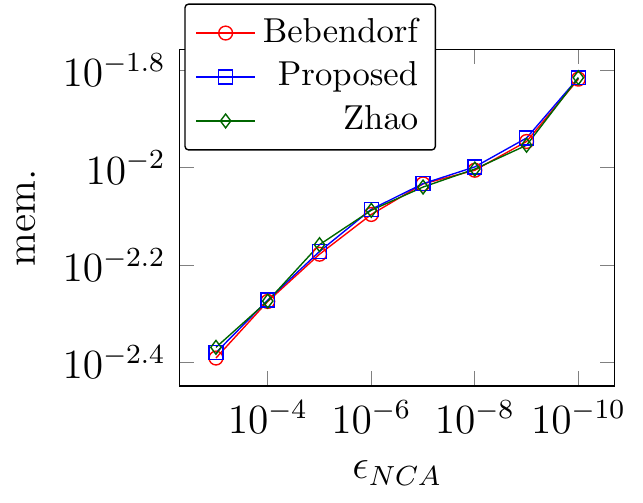}
        \end{subfigure}
        \hfill
                    \begin{subfigure}[b]{0.45\textwidth}
            \centering
            \includegraphics{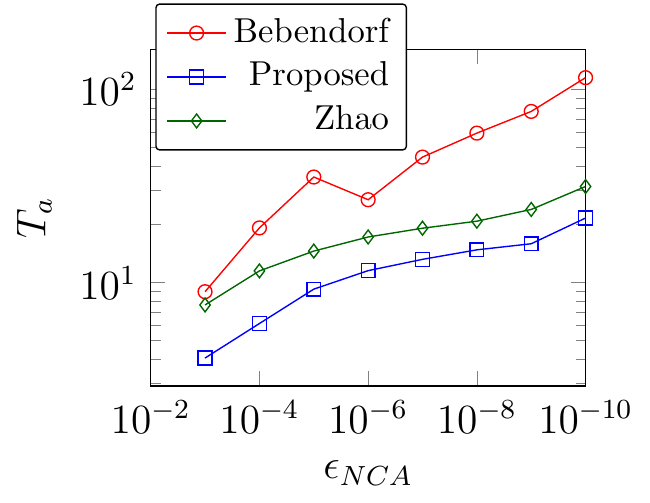}
            \end{subfigure}
            
        \begin{subfigure}[b]{0.45\textwidth}
          \centering
          \includegraphics{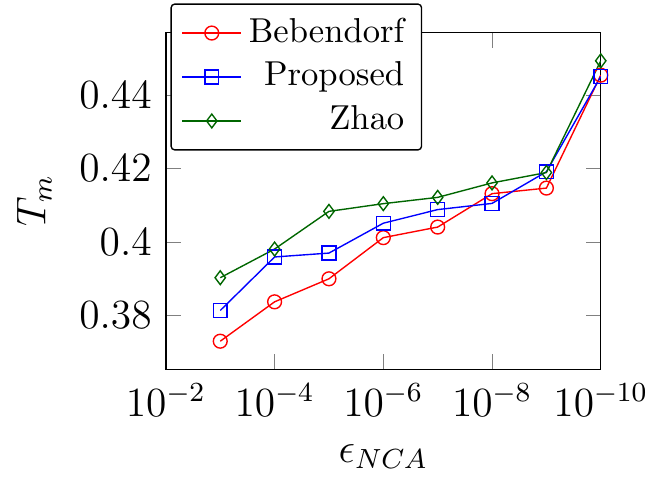}
          \end{subfigure}%
                  \hfill
          \begin{subfigure}[b]{0.45\textwidth}
            \centering
            \includegraphics{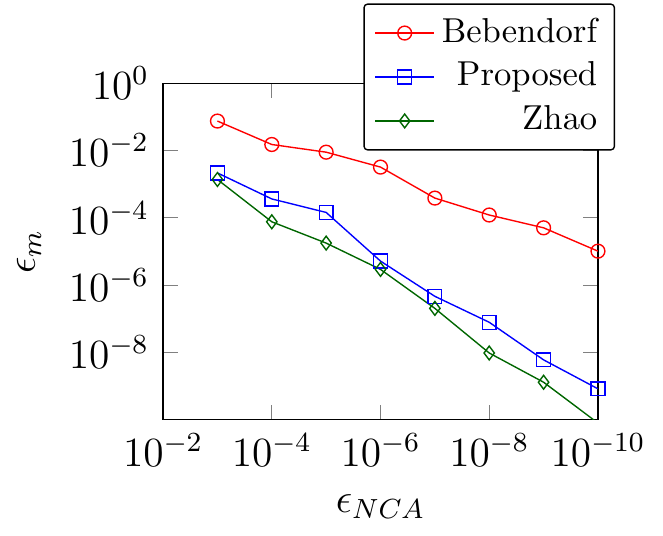}
            \end{subfigure}
            
      \caption{Results obtained with Experiment 1. Plots of memory, assembly time, matrix-vector product time, and relative error versus $\epsilon_{NCA}$ for kernel $\left(\frac{r(\log(r)-1)}{a(\log(a)-1)}\right)\chi_{r<a}+\left(\frac{\log(r)}{\log(a)}\right)\chi_{r\geq a}$ with $N=102400$.}
      \label{1VsEpsilon}
  \end{center}
\end{figure}

\begin{figure}[H]
  \begin{center}
      \begin{subfigure}[b]{0.45\textwidth}
        \centering
                    \includegraphics{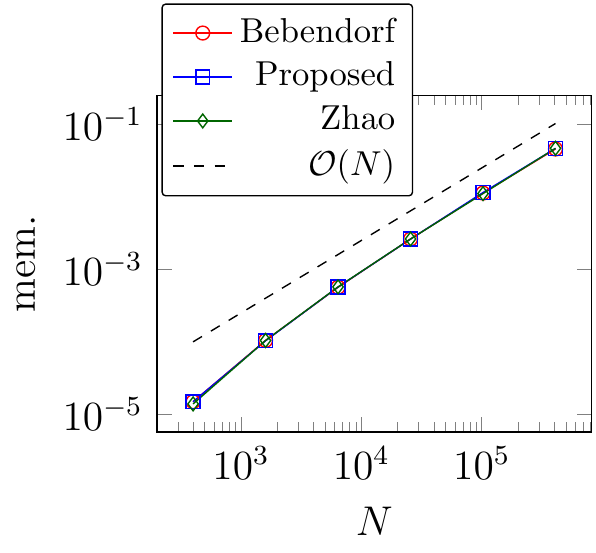}
        \end{subfigure}
                \hfill
      \begin{subfigure}[b]{0.45\textwidth}
        \centering
        \includegraphics{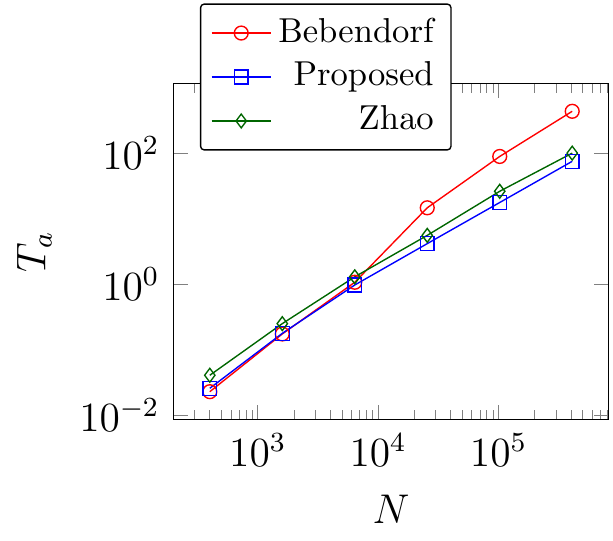}
        \end{subfigure}%
        
        \begin{subfigure}[b]{0.45\textwidth}
          \centering
          \includegraphics{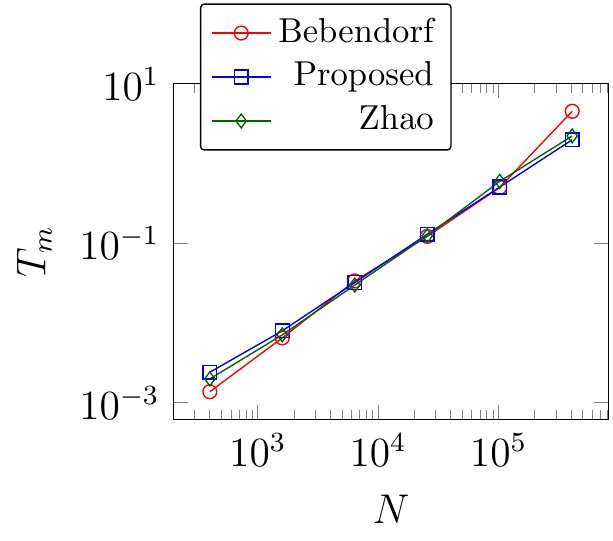}
          \end{subfigure}%
                  \hfill
    \begin{subfigure}[b]{0.45\textwidth}
            \centering
            \includegraphics{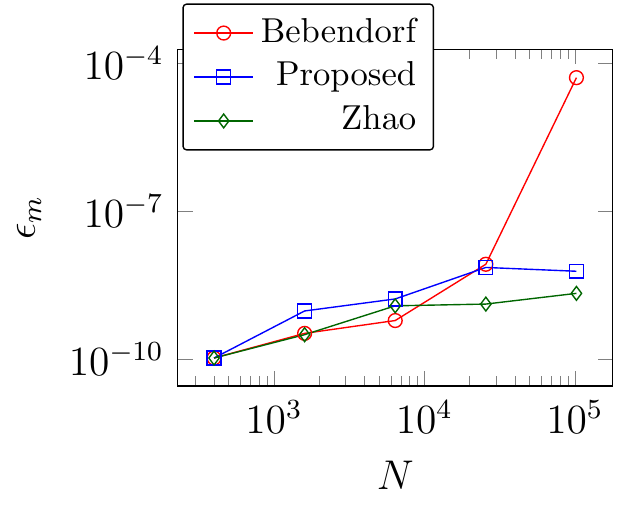}
            \end{subfigure}%
      \caption{Results obtained with Experiment 1. Plots of memory, assembly time, matrix-vector product time, and relative error versus $N$ for kernel $\left(\frac{r(\log(r)-1)}{a(\log(a)-1)}\right)\chi_{r<a}+\left(\frac{\log(r)}{\log(a)}\right)\chi_{r\geq a}$ with $\epsilon_{NCA}=10^{-9}$.}
      \label{1VsN}
  \end{center}
\end{figure}

\begin{figure}[H]
  \begin{center}
      \begin{subfigure}[b]{0.45\textwidth}
        \centering
                    \includegraphics{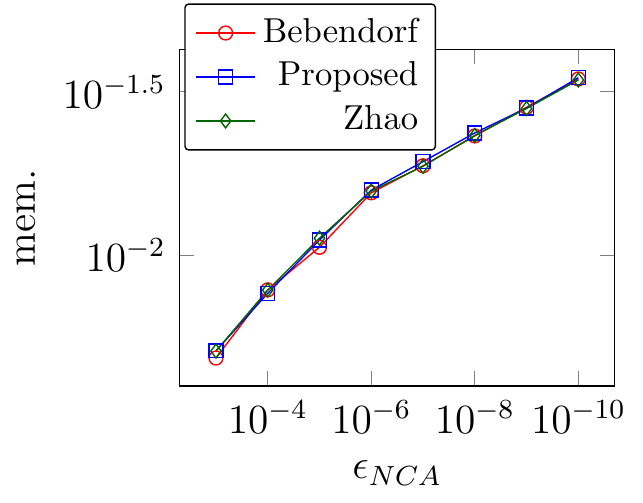}
        \end{subfigure}
                \hfill
        \begin{subfigure}[b]{0.45\textwidth}
        \centering
        \includegraphics{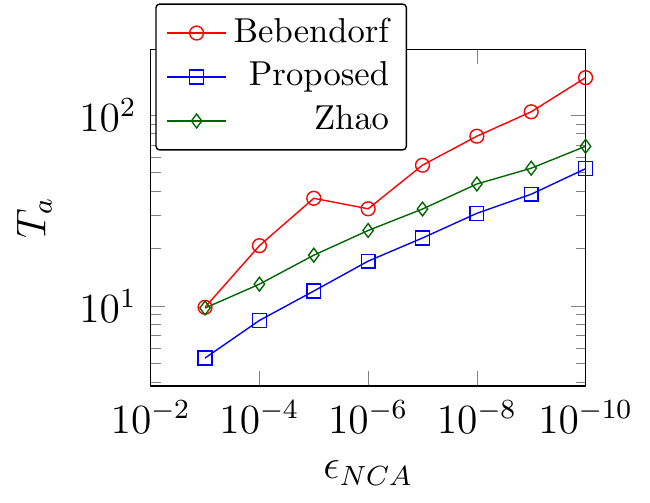}
        \end{subfigure}%
        
        \begin{subfigure}[b]{0.45\textwidth}
          \centering
          \includegraphics{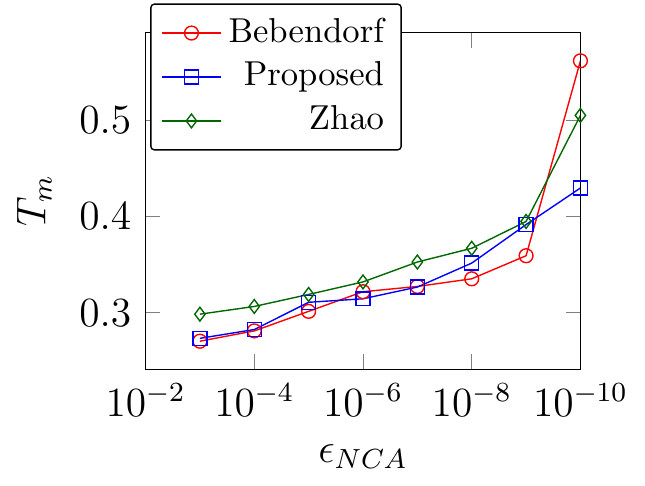}
          \end{subfigure}%
                  \hfill
          \begin{subfigure}[b]{0.45\textwidth}
            \centering
            \includegraphics{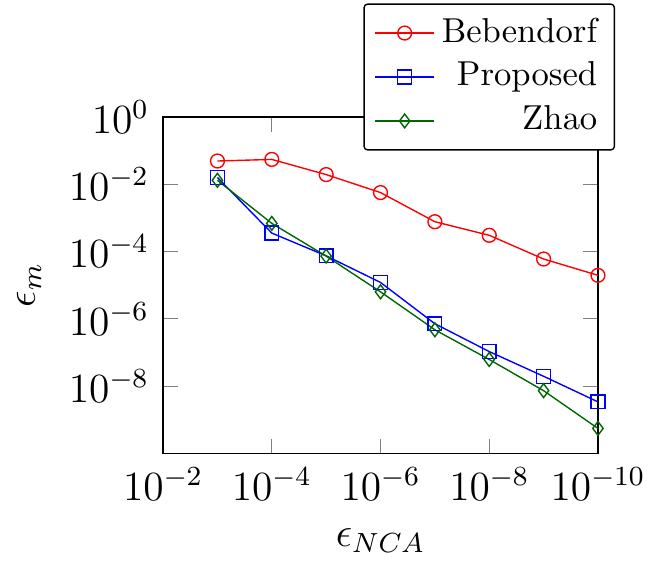}
            \end{subfigure}%
      \caption{Results obtained with Experiment 1. Plots of memory, assembly time, matrix-vector product time, and relative error versus $\epsilon_{NCA}$ for kernel $\left(\frac{r}{a}\right)\chi_{r<a}+\left(\frac{a}{r}\right)\chi_{r\geq a}$ with $N=102400$.}
      \label{2VsEpsilon}
  \end{center}
\end{figure}

\begin{figure}[H]
  \begin{center}
      \begin{subfigure}[b]{0.45\textwidth}
        \centering
                    \includegraphics{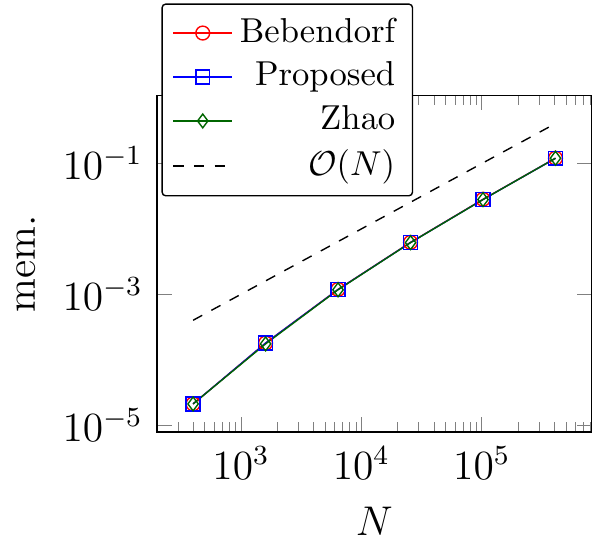}
        \end{subfigure}
                \hfill
      \begin{subfigure}[b]{0.45\textwidth}
        \centering
        \includegraphics{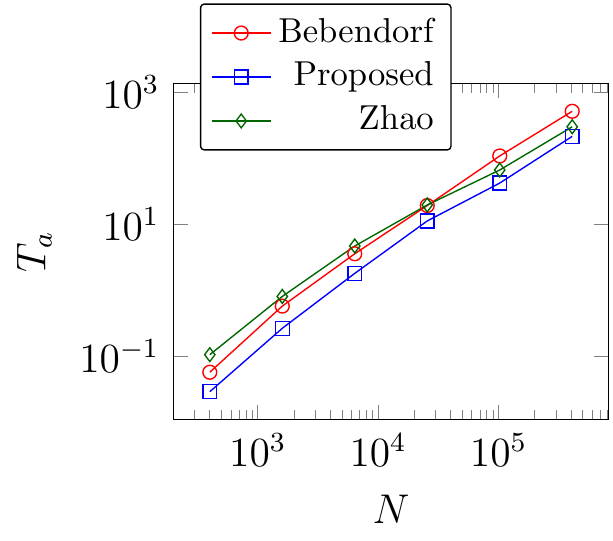}
        \end{subfigure}%
        
        \begin{subfigure}[b]{0.45\textwidth}
          \centering
          \includegraphics{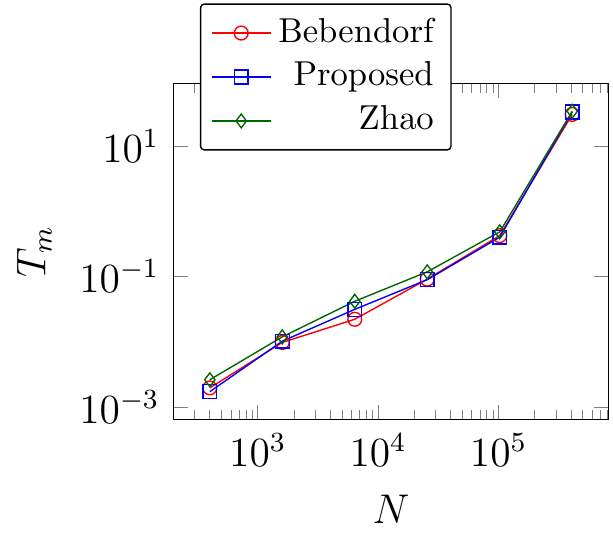}
          \end{subfigure}%
                  \hfill
          \begin{subfigure}[b]{0.45\textwidth}
            \centering
            \includegraphics{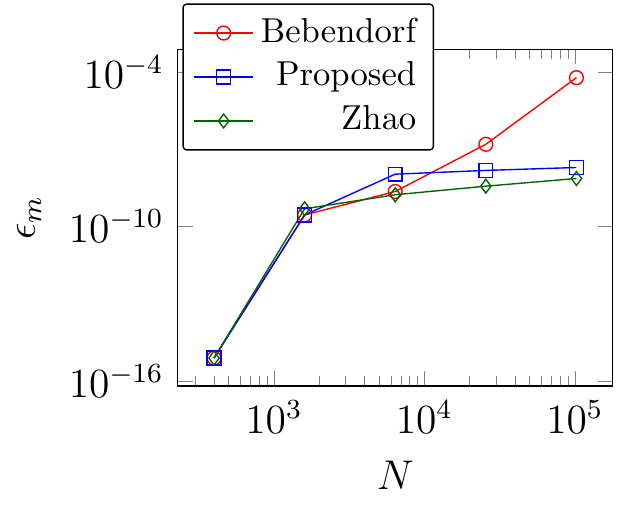}
            \end{subfigure}%
      \caption{Results obtained with Experiment 1. Plots of memory, assembly time, matrix-vector product time, and relative error versus $N$ for kernel $\left(\frac{r}{a}\right)\chi_{r<a}+\left(\frac{a}{r}\right)\chi_{r\geq a}$ with $\epsilon_{NCA}=10^{-9}$.}
      \label{2VsN}
  \end{center}
\end{figure}

Following observations are to be made from Figures~\ref{1VsEpsilon}, ~\ref{1VsN}, ~\ref{2VsEpsilon}, and ~\ref{2VsN}, which illustrate the scaling of mem., $T_{a}$, $T_{m}$, and $\epsilon_{m}$ with $\epsilon_{NCA}$ and $N$.
When scaling with $N$ is studied, $\epsilon_{NCA}$ is fixed to $10^{-9}$.
And when scaling with $\epsilon_{NCA}$ is studied, $N$ is fixed to $102400$.
\begin{enumerate}
    \item 
    mem. is almost same for the three algorithms.
    \item 
    The assembly time of the algorithm of~\cite{bebendorf2012constructing} scales as $\mathcal{O}(N\log(N))$. The assembly times of the NCA of~\cite{zhao2019fast} and NNCA scale as $\mathcal{O}(N)$. 
    \item
    NNCA is the fastest in terms of assembly. It is to be noted that this observation is made by ensuring that the relative error of NNCA either doesn't substantially differ from that of the existing NCAs or is lower than that of the existing NCAs. 
    \item 
    The time complexity of the matrix-vector products of all the three algorithms is $\mathcal{O}(N)$.
    \item 
    The relative error of NNCA is nearly equal to that of the NCA of Zhao et al. For large values of $N$, it is to be observed that the accuracies of NNCA and NCA of Zhao et al. are better than that of the NCA of Bebendorf et al. It is also to be observed from Figures~\ref{1VsEpsilon} and~\ref{2VsEpsilon} that $\epsilon_{m}$ of NNCA decreases as $\epsilon_{NCA}$ decreases, validating the convergence of NNCA.

\end{enumerate}

\subsection{Experiment 2: Matrix-vector product with a Chebyshev distribution of points in 2D}\label{sec:exp2}
\begin{figure}[H]
  \begin{center}
    \begin{subfigure}[b]{0.45\textwidth}
      \centering
      \includegraphics{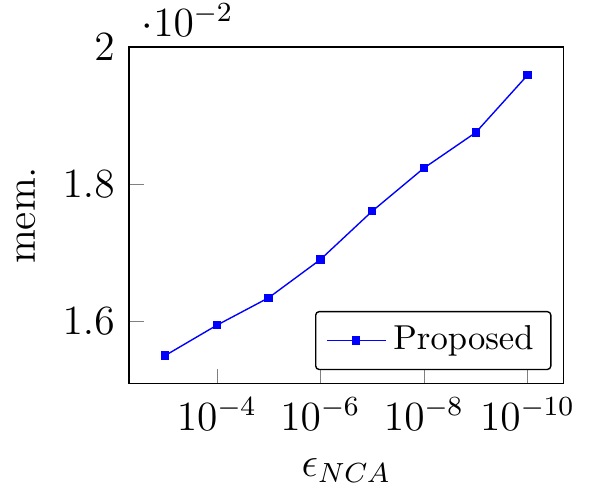}
      \end{subfigure}%
              \hfill
      \begin{subfigure}[b]{0.45\textwidth}
        \centering
        \includegraphics{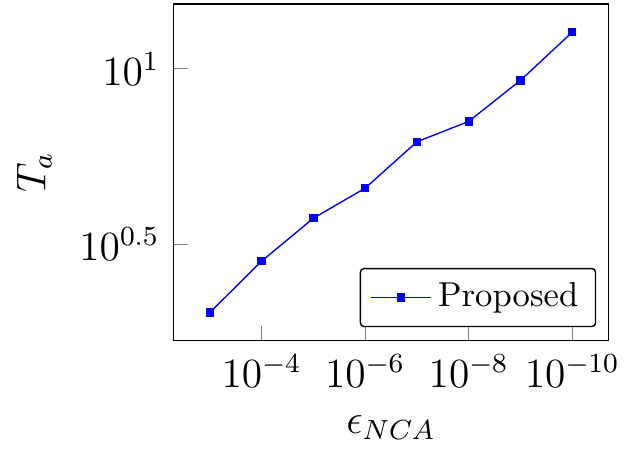}
        \end{subfigure}%

        \begin{subfigure}[b]{0.45\textwidth}
          \centering
           \includegraphics{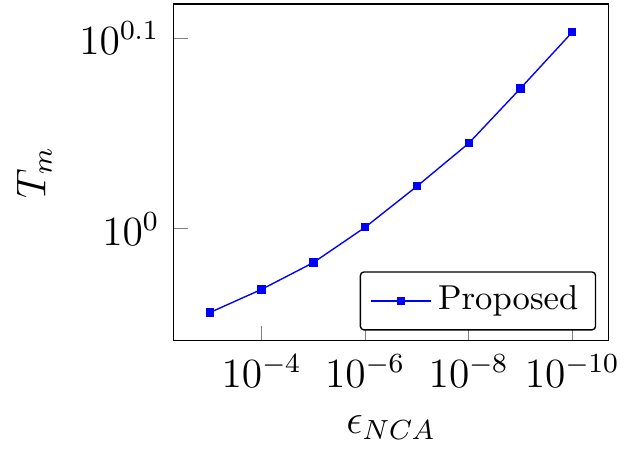}
          \end{subfigure}%
                  \hfill
          \begin{subfigure}[b]{0.45\textwidth}
            \centering
            \includegraphics{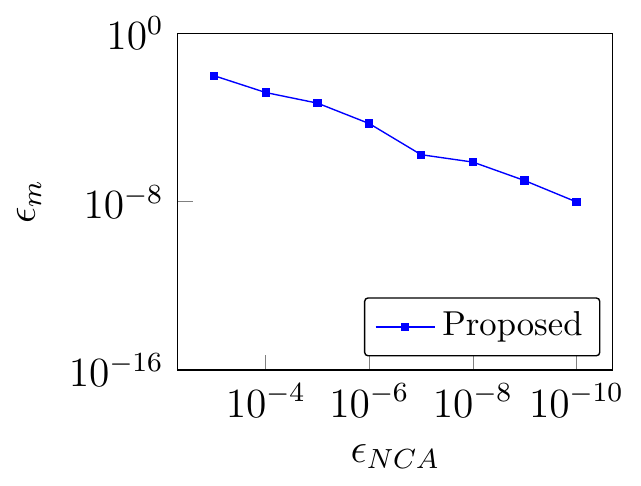}
            \end{subfigure}%
      \caption{Results obtained with Experiment 2. Plots of memory, assembly time, matrix-vector product time, and relative error versus $\epsilon_{NCA}$ for kernel $\left(\frac{r(\log(r)-1)}{a(\log(a)-1)}\right)\chi_{r<a}+\left(\frac{\log(r)}{\log(a)}\right)\chi_{r\geq a}$ with $N=102400$.}
      \label{nonUniformVsEps2D}
  \end{center}
\end{figure}

\begin{figure}[H]
  \begin{center}
    \begin{subfigure}[b]{0.45\textwidth}
      \centering
      \includegraphics{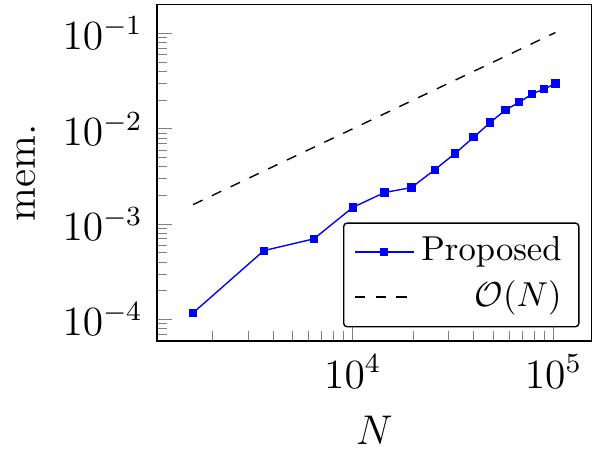}
      \end{subfigure}%
              \hfill
      \begin{subfigure}[b]{0.45\textwidth}
        \centering
        \includegraphics{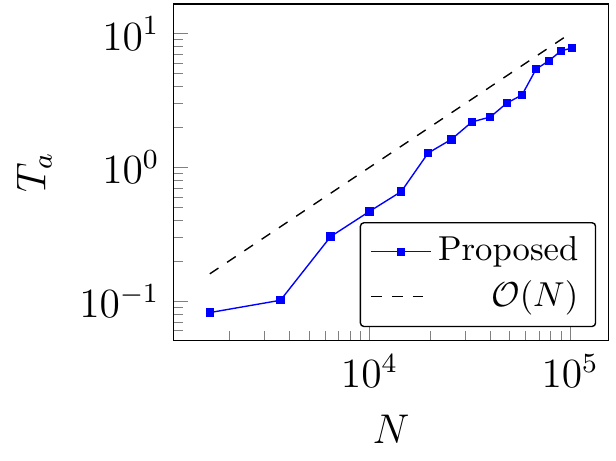}
        \end{subfigure}

        \begin{subfigure}[b]{0.45\textwidth}
          \centering
          \includegraphics{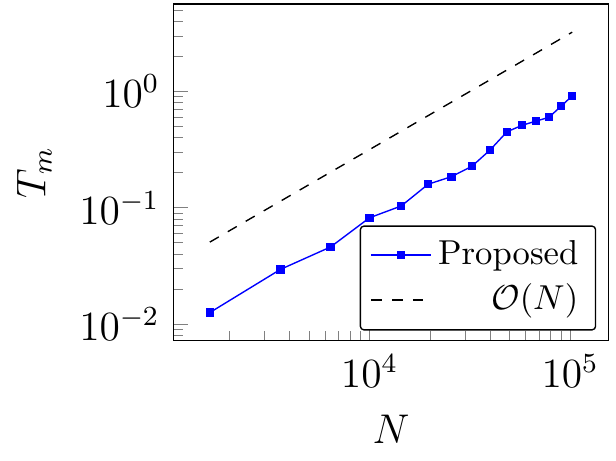}
          \end{subfigure}%
                  \hfill
          \begin{subfigure}[b]{0.45\textwidth}
            \centering
            \includegraphics{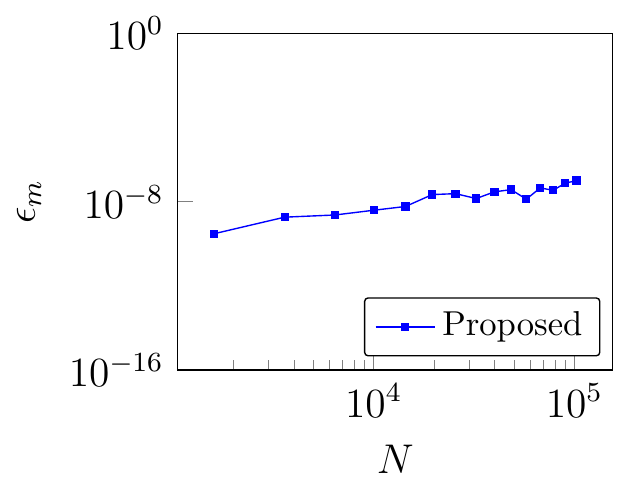}
            \end{subfigure}
      \caption{Results obtained with Experiment 2. Plots of memory, assembly time, matrix-vector product time, and relative error versus $N$ for kernel $\left(\frac{r(\log(r)-1)}{a(\log(a)-1)}\right)\chi_{r<a}+\left(\frac{\log(r)}{\log(a)}\right)\chi_{r\geq a}$ with $\epsilon_{NCA}=10^{-9}$.}
      \label{nonUniformVsN2D}
  \end{center}
\end{figure}

We perform this experiment to demonstrate the performance of NNCA based $\mathcal{H}^{2}$ matrix-vector in 2D when the particles follow a non-uniform distribution.
We consider the particles to be located on the tensor product Chebyshev grid of the domain $[-1,1]^{2}$.
Kernel function $\left(\frac{r(\log(r)-1)}{a(\log(a)-1)}\right)\chi_{r<a}+\left(\frac{\log(r)}{\log(a)}\right)\chi_{r\geq a}$ is considered.
In Figure~\ref{nonUniformVsEps2D}, we present various benchmarks as a function of $\epsilon_{NCA}$ with $N$ fixed to $102400$.
In Figure~\ref{nonUniformVsN2D}, we present various benchmarks as a function of $N$, with $\epsilon_{NCA}$ fixed to $10^{-9}$.
It is to be observed from Figure~\ref{nonUniformVsEps2D} that $\epsilon_{m}$ decreases as $\epsilon_{NCA}$ decreases, validating the convergence of NNCA.
It is also to be observed from Figure~\ref{nonUniformVsN2D} that mem., $T_{a}$, and $T_{m}$ scale linearly with $N$.

\subsection{Experiment 3: Matrix-vector product with a uniform distribution of points in 3D}
\begin{figure}[H]
  \begin{center}
    \begin{subfigure}[b]{0.45\textwidth}
      \centering
      \includegraphics{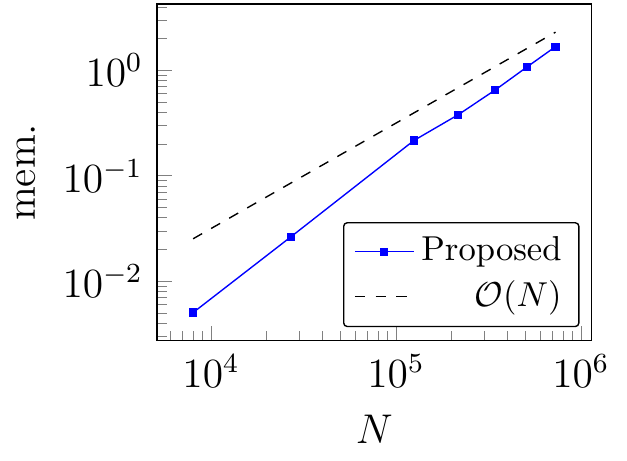}
      \end{subfigure}%
              \hfill
      \begin{subfigure}[b]{0.45\textwidth}
        \centering
        \includegraphics{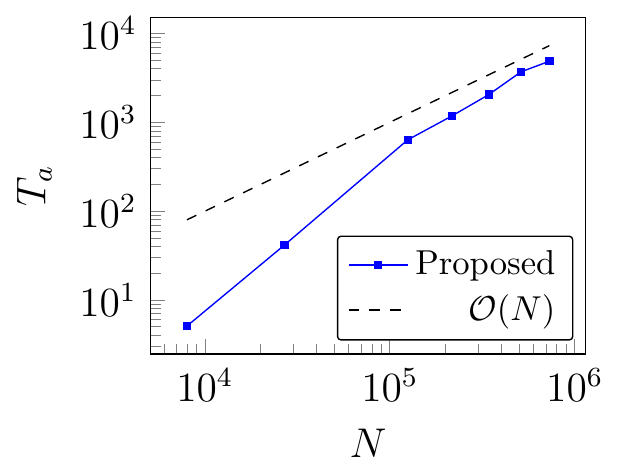}
        \end{subfigure}%

        \begin{subfigure}[b]{0.45\textwidth}
          \centering
          \includegraphics{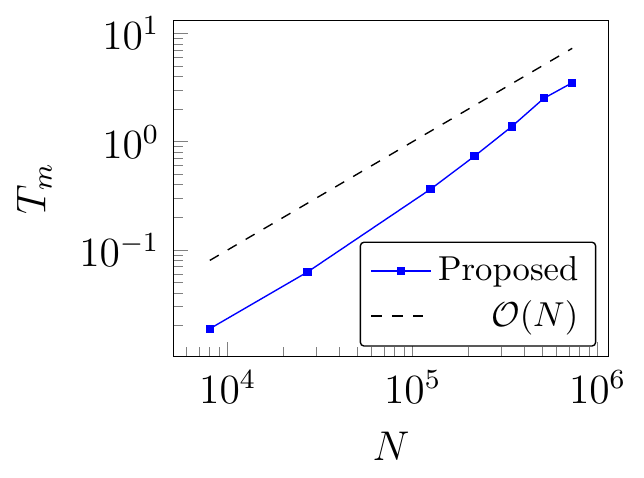}
          \end{subfigure}%
                  \hfill
          \begin{subfigure}[b]{0.45\textwidth}
            \centering
            \includegraphics{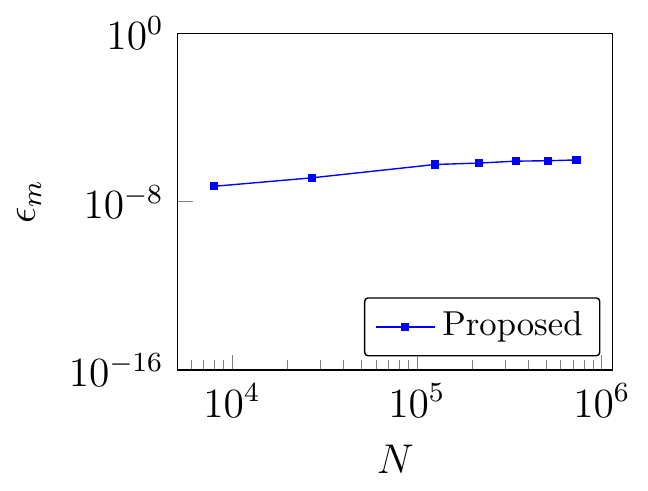}
            \end{subfigure}%
      \caption{Results obtained with Experiment 3. Plots of memory, assembly time, matrix-vector product time, and relative error versus $N$ for kernel function $\left(\frac{r}{a}\right)\chi_{r<a}+\left(\frac{a}{r}\right)\chi_{r\geq a}$.}
      \label{fig:3D_matvec1}
  \end{center}
\end{figure}

Here we demonstrate the performance of NNCA based $\mathcal{H}^{2}$ matrix-vector product, for the kernel function $\left(\frac{r}{a}\right)\chi_{r<a}+\left(\frac{a}{r}\right)\chi_{r\geq a}$ in 3D. $\epsilon_{NCA}$ is set to $10^{-7}$. A uniform distribution of particles is considered in the domain $[-1,1]^{2}$. In Figure~\ref{fig:3D_matvec1}, we benchmark it for various values of $N$.

\subsection{Experiment 4: Matrix-vector product with a Chebyshev distribution of points in 3D}
\begin{figure}[H]
  \begin{center}
    \begin{subfigure}[b]{0.45\textwidth}
      \centering
      \includegraphics{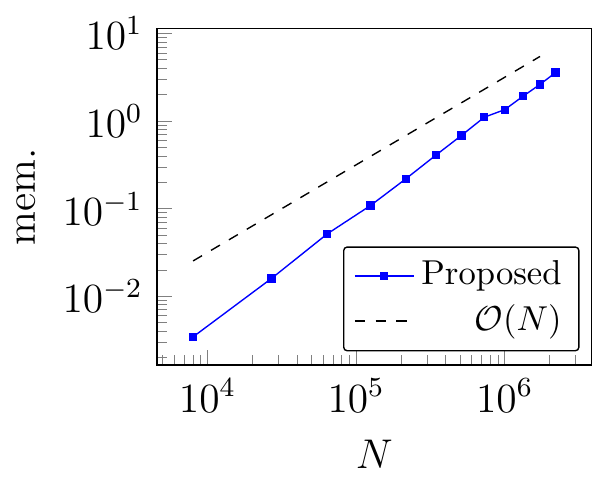}
      \end{subfigure}%
              \hfill
      \begin{subfigure}[b]{0.45\textwidth}
        \centering
        \includegraphics{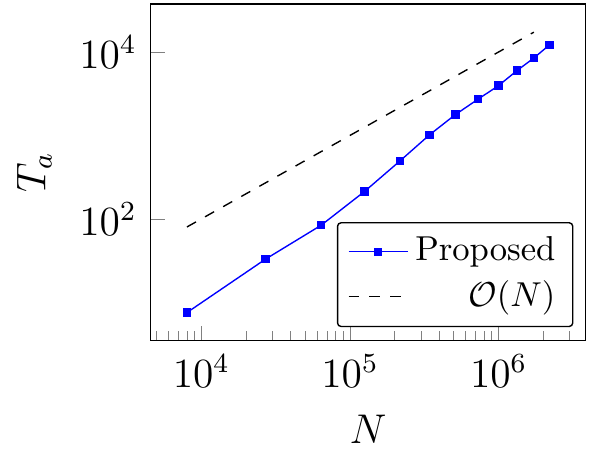}
        \end{subfigure}%

        \begin{subfigure}[b]{0.45\textwidth}
          \centering
          \includegraphics{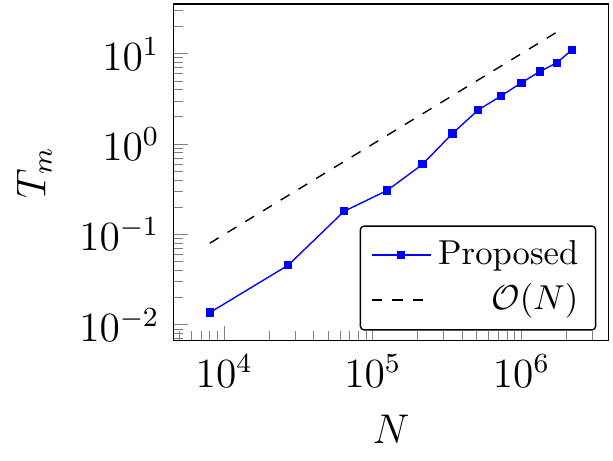}
          \end{subfigure}%
                  \hfill
          \begin{subfigure}[b]{0.45\textwidth}
            \centering
            \includegraphics{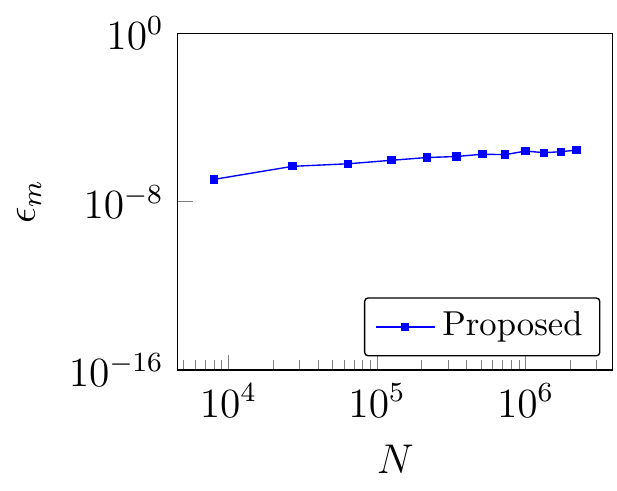}
            \end{subfigure}%
      \caption{Results obtained with Experiment 4. Plots of memory, assembly time, matrix-vector product time, and relative error versus $N$ for kernel function
     $\left(\frac{r}{a}\right)\chi_{r<a}+\left(\frac{a}{r}\right)\chi_{r\geq a}$.}
      \label{fig:exp5_matvec}
  \end{center}
\end{figure}

This experiment demonstrates the performance of NNCA based $\mathcal{H}^{2}$ matrix-vector product for the kernel function $\left(\frac{r}{a}\right)\chi_{r<a}+\left(\frac{a}{r}\right)\chi_{r\geq a}$. We consider the particles to be located on the tensor product Chebyshev grid of the domain $[-1,1]^{3}$. $\epsilon_{NCA}$ is set to $10^{-7}$. In Figure~\ref{fig:exp5_matvec}, we illustrate $T_{a}$, $T_{m}$ and $\epsilon_{m}$ versus $N$. 
It is to be observed that mem., $T_{a}$, and $T_{m}$ scale linearly with $N$.

\subsection{Experiment 5: Integral equation solver in 3D}
\begin{figure}[H]
  \begin{center}
    \begin{subfigure}[b]{0.45\textwidth}
      \centering
      \includegraphics{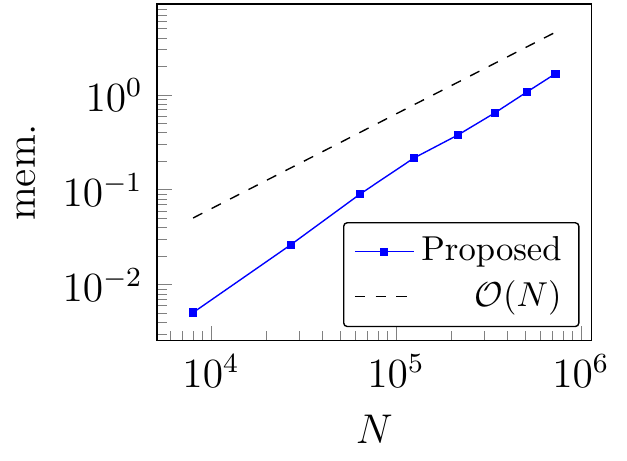}
      \end{subfigure}%
              \hfill
      \begin{subfigure}[b]{0.45\textwidth}
        \centering
        \includegraphics{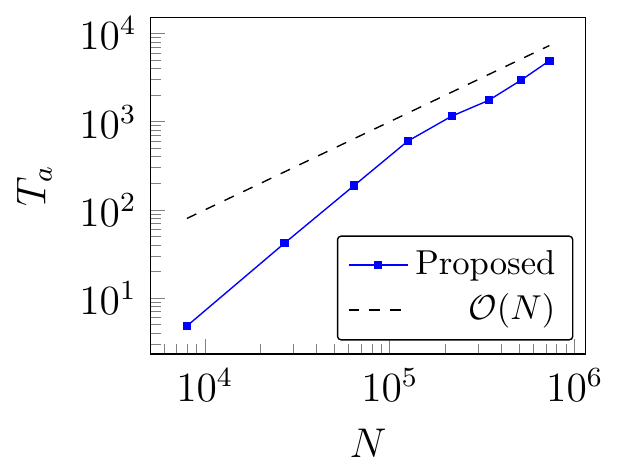}
        \end{subfigure}%

        \begin{subfigure}[b]{0.45\textwidth}
          \centering
          \includegraphics{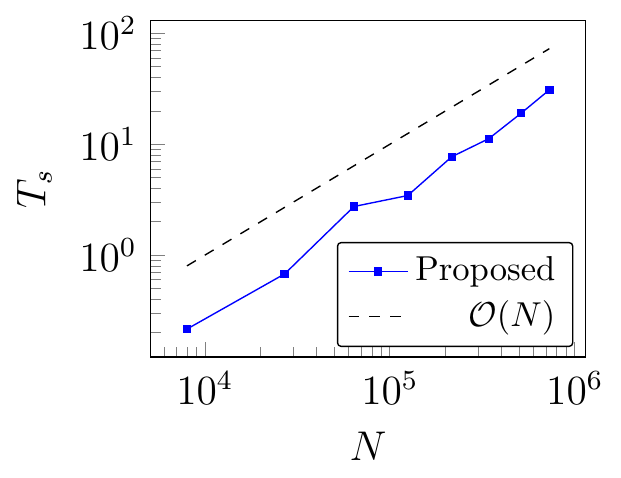}
          \end{subfigure}%
                  \hfill
          \begin{subfigure}[b]{0.45\textwidth}
            \centering
            \includegraphics{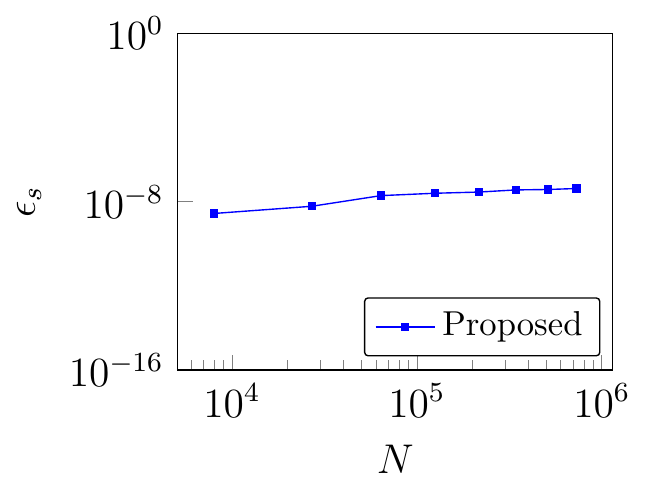}
            \end{subfigure}%
      \caption{Results obtained with Experiment 5. Plots of memory, assembly time, solve time, and relative error versus $N$.}
      \label{fig:3D_matvec2}
  \end{center}
\end{figure}

We consider the Fredholm integral equation of the second kind in the domain $[-1,1]^{3}$,
\begin{equation}
    \sigma(x) + \int_{D}K(x,y)\sigma(y)dy = f(x),\label{eq:Fredholm}
\end{equation}
with $K(x,y)=\frac{1}{\lVert x-y\rVert_{2}}$. We consider a Nystrom discretization of the integral equation~\ref{eq:Fredholm}, on a uniform grid, to obtain a discrete linear system of the form 
\begin{equation}
A\vec{\sigma}=\vec{f}.\label{eq:discreteFredholm}
\end{equation}
We consider a random vector $\vec{\sigma}$ and find the vector $\vec{f}=A\vec{\sigma}$ (exact upto roundoff). 
With $\vec{f}$ as the right hand side in equation~\ref{eq:discreteFredholm}, 
we solve for $\vec{\sigma}$ using GMRES, where the matrix-vector product to be performed in each iteration of GMRES is computed using NNCA based $\mathcal{H}^{2}$ matrix-vector product. Let $\vec{\sigma_{t}}$ be the computed $\vec{\sigma}$. $\epsilon_{NCA}$ is set to $10^{-7}$. $\epsilon_{GMRES}$, the stopping criterion for GMRES, is set to $10^{-10}$, i.e., we stop further iterations if the relative residual is less than $\epsilon_{GMRES}$. We define the relative forward error, $\epsilon_{s}$, as $\frac{||\vec{\sigma_{t}}-\vec{\sigma}||_{2}}{||\vec{\sigma}||_{2}}$. 
In Table~\ref{table:exp6pVsN} and Figure~\ref{fig:3D_matvec2}, we illustrate iter., mem., $T_{a}$, $T_{s}$, and $\epsilon_{s}$ for various values of $N$.

\begin{table}[H]
\centering
  \begin{tabular}{|c|c|c|c|c|c|c|c|c|c|c|c|c|}
    \hline
    $N$ &&&&&&&&\\
    (in thousands) & 8 & 27 & 64 & 125 & 216 & 343 & 512 & 729\\ \hline \hline
    iter. & 8 & 8 & 8 & 8 & 7 & 7 & 7 & 7\\ \hline
  \end{tabular}
 \caption{iter. versus $N$ obtained with Experiment 5.}
    \label{table:exp6pVsN}
 \end{table}

\subsection{Experiment 6: Kernel SVM in 4D}
Support Vector Machine (SVM)~\cite{boser1992training} is a well-known method belonging to the class of supervised machine learning methods. It is widely used as a classification algorithm. The naive SVM can classify only linearly separable data. In order to classify data that is not linearly separable, kernel SVM is used. 

Let $\{x_{i}\}_{i=1}^{N}$ be $N$ data points available to train the model belonging to one of the two classes identified as class 1 and class 2. Let $y_{i}$ represent the label associated with data point $x_{i}$, defined as
    \[y_i = \begin{cases}
    1 & \text{ if } x_{i} \in \text{class 1} \\
    -1 & \text{ if } x_{i} \in \text{class 2} 
    \end{cases}.\]
Kernel SVM involves maximizing the following objective function
\begin{equation} \label{svm_eq1}
    \underset{\alpha}{\operatorname{argmax}} \dsum_{i=1}^N \alpha_i - \dfrac{1}{2} \dsum_{i=1}^N \dsum_{j=1}^N K \bkt{x_i, x_j} y_i y_j \alpha_i \alpha_j
\end{equation}
subject to $0\leq\alpha\leq\lambda$ and $\dsum_{i=1}^{N}\alpha_{i}y_{i}=0$. Here $\alpha=[\alpha_{1}, \alpha_{2}, \hdots, \alpha_{N}]$ and $K(x_{i}, x_{j})$ is the kernel function evaluation at points $x_{i}$ and $x_{j}$. Some of the widely used kernel functions are Gaussian kernel, Laplace kernel, sigmoid kernel, polynomial kernel, Matérn kernel, etc. 

The optimization problem in Equation~\eqref{svm_eq1} can be solved using Lagrange multipliers. Let $L(\alpha)$ be defined as
\begin{equation} \label{svm_eq2}
    L\bkt{\alpha} = \dsum_{i=1}^N \alpha_i - \dfrac{1}{2} \dsum_{i=1}^N \dsum_{j=1}^N K \bkt{x_i, x_j} y_i y_j \alpha_i \alpha_j - \dfrac{1}{2} \beta \bkt{\dsum_{j=1}^N \alpha_j y_j}^2
\end{equation}
The $\alpha$ that maximizes $L(\alpha)$ can be found using gradient descent. Let $\eta$ be the learning rate. $\alpha_{i}$ is found iteratively as follows
\begin{equation} \label{svm_eqgd}
    \alpha_i := \alpha_i + \eta \dfrac{\partial L\bkt{\alpha}}{\partial \alpha_i}
\end{equation}
where
\begin{equation} \label{svm_eq3}
    \dfrac{\partial L\bkt{\alpha}}{\partial \alpha_i} = 1 - \dsum_{j=1}^N K \bkt{x_i, x_j} y_i y_j \alpha_j -  \beta \dsum_{j=1}^N \alpha_j y_j y_i
\end{equation}
Equation~\eqref{svm_eq3} in matrix-vector parlance, expressed using MATLAB's notations, takes the following form
\begin{equation} \label{svm_eq5}
    DL \bkt{\alpha} = 1 - y \cdot* \bkt{K v} - \beta * \text{sum}(v) * y
\end{equation}
where
$DL \bkt{\alpha} = 
    \begin{bmatrix}
        \dfrac{\partial L\bkt{\alpha}}{\partial \alpha_1} \\
        \vdots  \\
        \dfrac{\partial L\bkt{\alpha}}{\partial \alpha_N} 
    \end{bmatrix},
    y = 
    \begin{bmatrix}
        y_1  \\
        \vdots  \\
        y_N
    \end{bmatrix},
    v=
    \begin{bmatrix}
        y_1 \alpha_1 \\
        \vdots  \\
        y_N \alpha_N
    \end{bmatrix}
$
and
$K = 
    \begin{bmatrix}
        K(x_1,x_1) & \dots  & K(x_1,x_N) \\
        \vdots & \vdots & \vdots \\
       K(x_N,x_1) & \dots  & K(x_N,x_N) \\
    \end{bmatrix}.
$
Equation~\eqref{svm_eqgd} is iterated multiple times until convergence. It involves the computation of a matrix-vector product $Kv$, in each iteration. This is computationally expensive when done naively. To reduce the computational time and memory, we propose to use the NNCA based $\mathcal{H}^{2}$ matrix-vector product. We term the SVM coupled with the NNCA based $\mathcal{H}^{2}$ matrix-vector product to be the Fast SVM (FSVM).

We implemented FSVM in C++ in $d$ dimensions. We illustrate the numerical benchmarks that we observed on synthetic datasets in 2D and 4D. We considered $M$ randomly distributed particles to be the feature vectors, belonging to $[-1.4,1.4]^{2}$ in 2D and $[-1.0,1.0]^{4}$ in 4D. Each of the feature vectors is assigned either class $1$ or $2$ such that the number of feature vectors belonging to classes 1 and 2 are nearly equal. Further, the dataset is divided into two parts: train data and test data. $85\%$ of the data is considered to be the train data and the rest is considered to be test data. It is also ensured that the number of data points belonging to classes 1 and 2 are nearly equal in both the train and test data. 

We identify the label of the test data after training. For $x$ belonging to the test data, we find $f(x) = \dsum_{i=1}^{N} \alpha_i y_i K \bkt{x_i, x} + b$, where $b$ is the bias. We then find the label of $x$ to be $g(x)= \text{sgn}\bkt{f(x)}$. We introduce some notations in Table~\ref{table:SVMNotations} to describe the numerical benchmarks of FSVM. 
\begin{table}[H]
  \centering
  \begin{tabular}{|l|p{0.8\textwidth}|}
    \hline
    $M$ & Total number of data points including train and test data \\
    \hline
    $C_{i}$ & number of train data points belonging to class $i$, where $i\in \{1,2\}$ \\
    \hline
    $c_{i}$ & number of test data points belonging to class $i$, where $i\in \{1,2\}$ \\
    \hline
    $s_{i}$ & number of test data points identified by SVM to be belonging to class $i$, where $i\in \{1,2\}$ \\
    \hline
    $t_{F}$ & Training time taken by FSVM that includes the time to assemble the Kernel matrix using NNCA\\
 \hline
 $t_{N}$ & Training time taken by normal SVM (NSVM)\\
\hline
  iter & Number of iterations taken by gradient descent for convergence.\\
 \hline
$i_{F}$ & Training time taken by FSVM per iteration $i_{F} = t_{F}/iter$ \\
 \hline
$i_{N}$ & Training time taken by NSVM per iteration $i_{N} = t_{N}/iter$ \\
  \hline
$A_{1}$ & Accuracy of class 1 test data $A_{1} = \frac{s_{1}}{c_{1}}\times 100\%$ \\
  \hline
$A_{2}$ & Accuracy of class 2 test data $A_{2} = \frac{s_{2}}{c_{2}}\times 100\%$ \\
   \hline
$OA$ & Accuracy of test data $A_{1} = \frac{s_{1}+c_{1}}{s_{2}+c_{2}}\times 100\%$ \\
\hline
\end{tabular}
\caption{List of notations followed in this subsection}
   \label{table:SVMNotations}
\end{table}

\begin{figure}[H]
  \begin{center}
      \begin{subfigure}[b]{0.45\textwidth}
        \centering
    \includegraphics[scale=0.25]{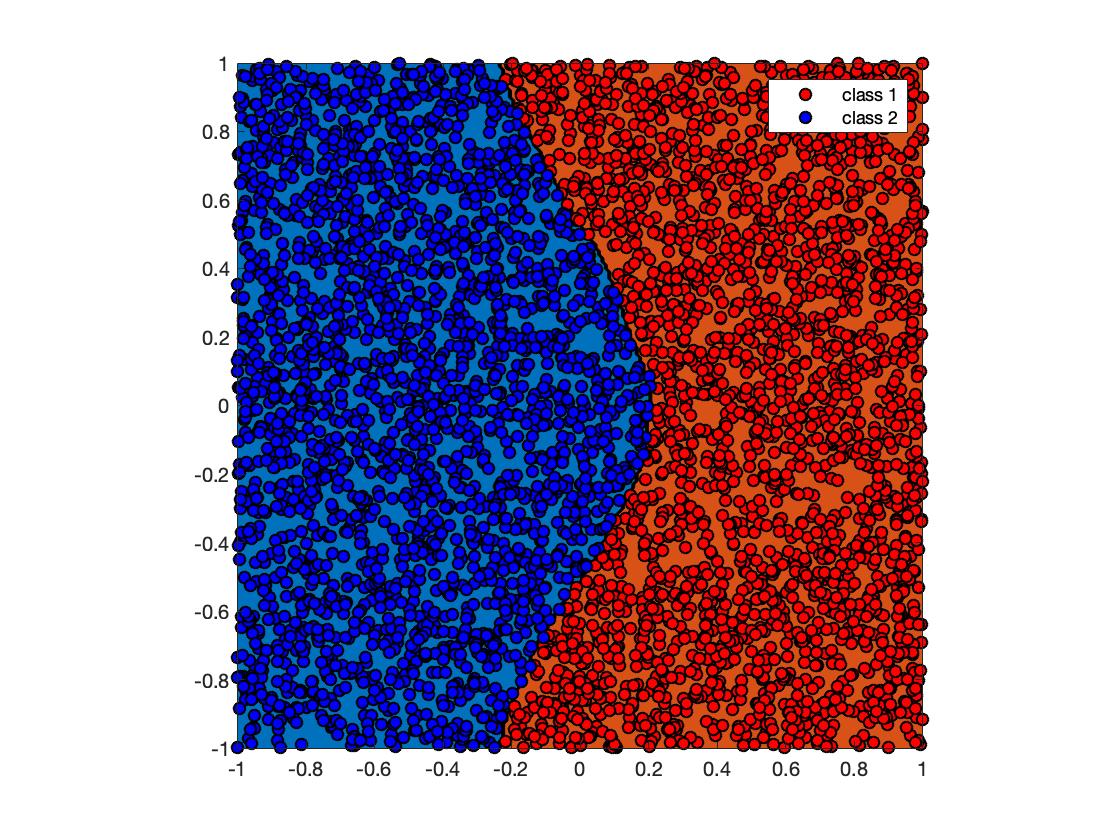}
        \caption{}
        \label{fig:svm1}
        \end{subfigure}
        \hfill
                    \begin{subfigure}[b]{0.45\textwidth}
            \centering
            \includegraphics[scale=0.25]{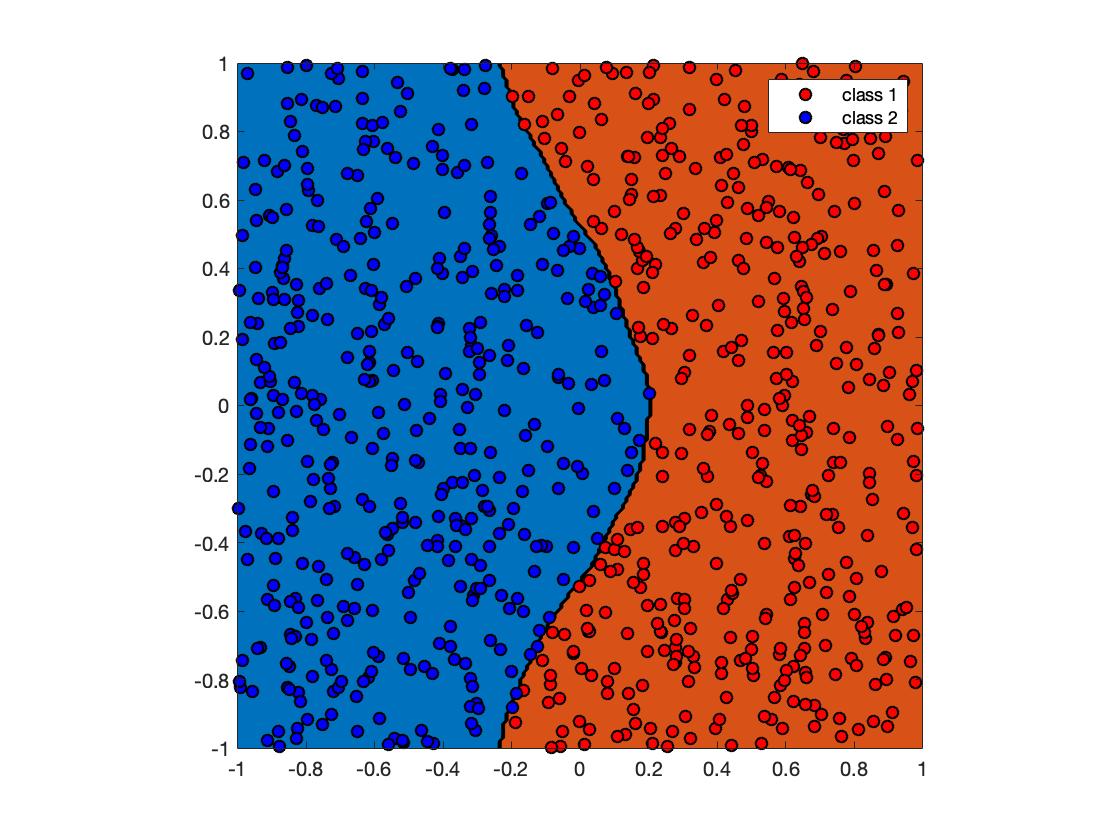}
            \caption{}
            \label{fig:svm2}
            \end{subfigure}
            
        \begin{subfigure}[b]{0.45\textwidth}
          \centering
          \includegraphics[scale=0.25]{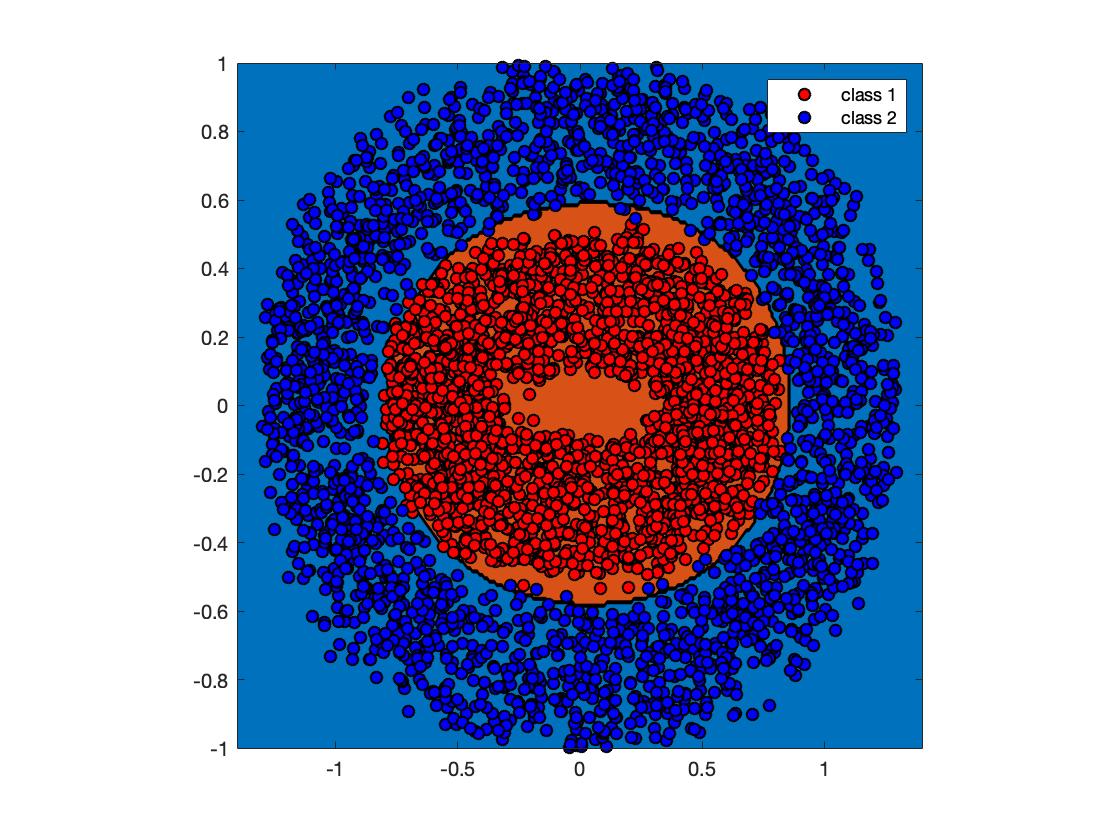}
                 \caption{}
            \label{fig:svm3}
                 \end{subfigure}%
                  \hfill
          \begin{subfigure}[b]{0.45\textwidth}
            \centering
            \includegraphics[scale=0.25]{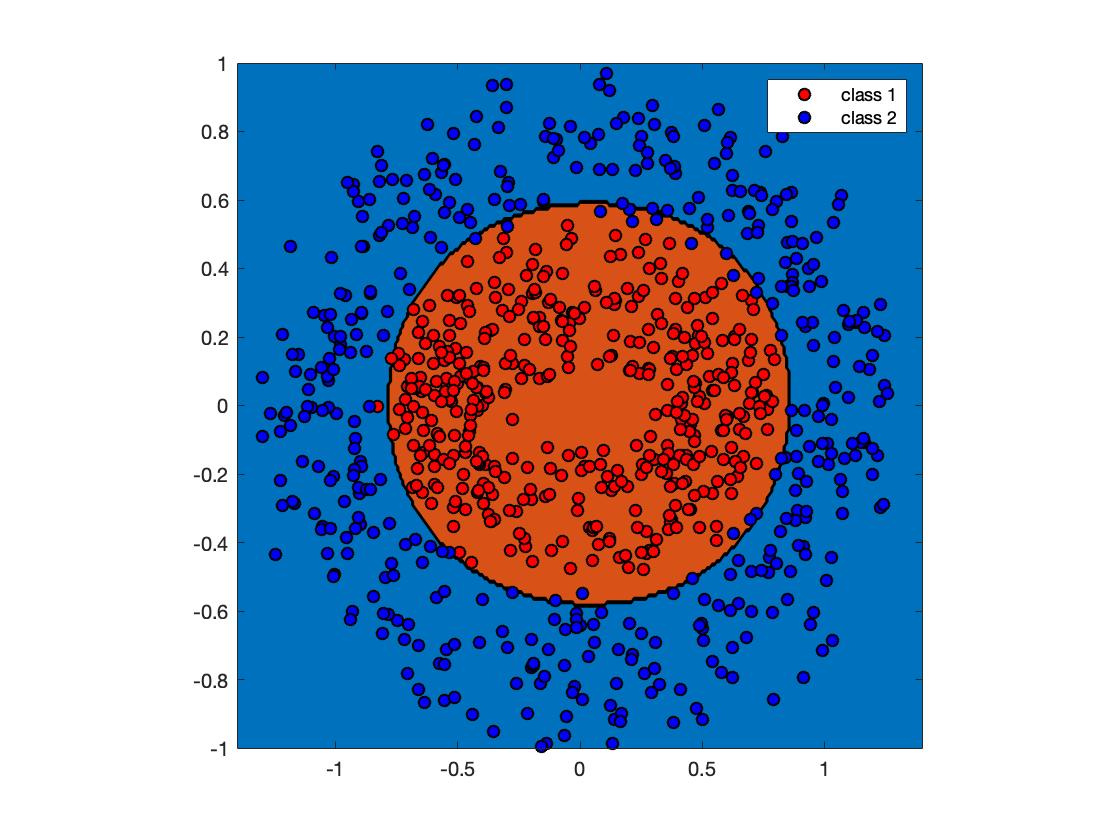}
            \caption{}
            \label{fig:svm4}
            \end{subfigure}
            
      \caption{Results obtained with Experiment 6 in 2D. (a) and (b) Decision boundary constructed using FSVM with Matérn kernel superposed on the train and test datasets of Dataset 1 respectively. (c) and (d) Decision boundary constructed using FSVM with Gaussian kernel superposed on the train and test datasets of Dataset 2 respectively.}
      \label{svm:decision2D}
  \end{center}
\end{figure}

NSVM is the normal SVM, where the matrix-vector products are computed naively. In Table~\ref{tab:SVM_2d}, we compared the performance of FSVM with the normal SVM (NSVM) on synthetic datasets with two features. We constructed two synthetic datasets Dataset 1 and Dataset 2. For Dataset 1, we used Matérn kernel, i.e. $K(x,y) = \exp(-\|x-y\|_{2})$ and for Dataset 2, we used Gaussian kernel, i.e. $K(x,y) = \exp(-\|x-y\|^{2}_{2})$. The speed-up of FSVM over NSVM can be observed from Table~\ref{tab:SVM_2d}. In Figure~\ref{svm:decision2D}, we illustrate the decision boundary for the two datasets.

In Table~\ref{tab:SVM_4}, we compared the performance of FSVM with the normal SVM (SVM where the matrix-vector products are computed naively, which we refer to as NSVM) on synthetic datasets with four features. We considered the Matérn kernel. The speed-up of FSVM over NSVM can be observed from the tables.

 \begin{table}[H]
 \begin{center}
 \resizebox{\textwidth}{!}{  
\begin{tabular}{|l|l|l|ll|ll|ll|ll|l|l|l|}
\hline
Kernel & Dataset & $M$ & \multicolumn{2}{c|}{Train data} & \multicolumn{2}{l|}{Test data} & \multicolumn{4}{c|}{Train time}    & \multirow{2}{*}{A$_1$} & \multirow{2}{*}{A$_2$} & \multirow{2}{*}{OA} \\ \cline{1-11}
& & & \multicolumn{1}{l|}{$C_1$}      & $C_2$     & \multicolumn{1}{l|}{$c_1$}    & $c_2$   & \multicolumn{1}{l|}{$t_{N}$} & $t_{F}$ & \multicolumn{1}{l|}{${i}_{N}$}  & ${i}_{F}$  &                      &                      &                    \\ \hline
Matérn & Dataset 1 &  5625 & \multicolumn{1}{l|}{2415}         &  2367       & \multicolumn{1}{l|}{426}       &  417     & \multicolumn{1}{l|}{702.2 }       & 27.4        & \multicolumn{1}{l|}{4.7}        &  0.18        &    100                  &       99.0               &      99.5              \\ \hline
Gaussian & Dataset 2 &  5625 & \multicolumn{1}{l|}{2381}         &  2400       & \multicolumn{1}{l|}{420}       &  423     & \multicolumn{1}{l|}{459.6}       & 20.0        & \multicolumn{1}{l|}{5.0}        &  0.22        &    99.5                  &       95.9               &      97.7              \\ \hline
\end{tabular}
}
\end{center}
     \caption{Results obtained with experiment 6 in 2D. Comparison of the performance of NSVM with FSVM using Gaussian kernel on a synthetic dataset.}
     \label{tab:SVM_2d}
\end{table}

 \begin{table}[H]
\begin{tabular}{|l|ll|ll|ll|ll|l|l|l|}
\hline
$M$ & \multicolumn{2}{c|}{Train data} & \multicolumn{2}{l|}{Test data} & \multicolumn{4}{c|}{Train time}    & \multirow{2}{*}{A$_1$} & \multirow{2}{*}{A$_2$} & \multirow{2}{*}{OA} \\ \cline{1-9}
  & \multicolumn{1}{l|}{$C_1$}      & $C_2$     & \multicolumn{1}{l|}{$c_1$}    & $c_2$   & \multicolumn{1}{l|}{$t_{N}$} & $t_{F}$ & \multicolumn{1}{l|}{${i}_{N}$}  & ${i}_{F}$  &                      &                      &                    \\ \hline
 4096 & \multicolumn{1}{l|}{1782}         &  1701       & \multicolumn{1}{l|}{314}       &  300     & \multicolumn{1}{l|}{1115.8}       & 203.0        & \multicolumn{1}{l|}{2.54}        &  0.46        &    99.4                  &       98.3               &      98.8              \\ \hline
 6561 & \multicolumn{1}{l|}{2811}         &  2767       & \multicolumn{1}{l|}{495}       &   488    & \multicolumn{1}{l|}{3662.9}       &  357.7       & \multicolumn{1}{l|}{6.25}        &  0.61        &   98.88                  &        98.9              &    98.8                \\ \hline
 10000 & \multicolumn{1}{l|}{4283}         &  4218       & \multicolumn{1}{l|}{755}       &  744     & \multicolumn{1}{l|}{11909.2}       &  694.1       & \multicolumn{1}{l|}{14.90}        &  0.86        &  99.2                    &     98.1                 &     98.6               \\ \hline
 14641 & \multicolumn{1}{l|}{6331}         &  6115       & \multicolumn{1}{l|}{1117}       &  1078     & \multicolumn{1}{l|}{30895.1}       &    1323.8     & \multicolumn{1}{l|}{31.78}        &  1.36        &  99.4                    &     98.2                 &     98.8               \\ \hline
\end{tabular}
     \caption{Results obtained with experiment 6 in 4D. Comparison of the performance of NSVM with FSVM using Matérn kernel on a synthetic dataset.}
     \label{tab:SVM_4}
\end{table}

\subsection{NNCA for oscillatory kernels}
In this article, NNCA is demonstrated for non-oscillatory kernels. We have also demonstrated in~\cite{gujjula2022new}, that our NNCA, with a modified admissibility condition, works for oscillatory kernels as well.
The 2D and 3D Helmholtz kernels follow a directional admissibility condition for low rank~\cite{engquist2009fast,engquist2010fast}. 
To develop NCA for the respective oscillatory kernels, one needs to adapt the admissibility condition to those of the oscillatory kernels.
For more details on the construction and benchmarks of NNCA for the 2D Helmholtz kernel, we refer the readers to~\cite{gujjula2022new}.

\section{Conclusion}
We proposed a new Nested Cross Approximation for $\mathcal{H}^{2}$ matrices and demonstrated its applicability by developing the $\mathcal{H}^{2}$ matrix-vector product.
The key highlight of NNCA is to choose the far-field pivots of a cell from its interaction list region instead of the entire far-field region.
We compared NNCA with the existing NCAs and demonstrated that NNCA outperforms the existing NCAs in the assembly time.
We further demonstrated the linear complexity of NNCA based $\mathcal{H}^{2}$ matrix-vector product by considering a comprehensive set of experiments in 2D and 3D. In addition, using the NNCA based $\mathcal{H}^{2}$ matrix-vector product, we accelerate i) solving an integral equation in 3D and ii) SVM on datasets with two and four features.
In the spirit of reproducible computational science, the implementation of the algorithm developed in this article is made available at \url{https://github.com/SAFRAN-LAB/NNCA}.

\section*{Acknowledgements}
The authors acknowledge HPCE, IIT Madras, India for providing access to the AQUA cluster. Vaishnavi Gujjula acknowledges the support of Women Leading IITM in Mathematics, IITM, India (SB22230053MAIITM008880). Sivaram Ambikasaran acknowledges the support of YSRA from BRNS, DAE, India (No.34/20/03/2017-BRNS/34278) and MATRICS grant from SERB, India (Sanction number: MTR/2019/001241).

\bibliographystyle{ACM-Reference-Format}
\bibliography{bibFile}

\end{document}

%% file: includes.tex
\usepackage{graphicx}
\usepackage{tikz}
\usepackage{subcaption}
\usepackage{float}
\usepackage{multirow}

\newcommand{\dsum}{\displaystyle\sum}
\newcommand{\bkt}[1]{\left(#1\right)}
\newtheorem{remark}{Remark}[section]